\def\no{\if01}
\def\iftwelvept{\no}

\def\ifusepdf{\no}
\def\ifpsfont{\no}

\iftwelvept
\documentclass[leqno,12pt]{amsart}
\else
\documentclass[leqno,10pt]{amsart}
\fi
\usepackage{amssymb}
\usepackage{amscd}
\usepackage{latexsym}
\usepackage{verbatim}
\usepackage[all]{xy}

\ifusepdf
\usepackage{hyperref}
\else\fi
\ifpsfont
\usepackage[T1]{fontenc}
\usepackage{times}
\else\fi

\iftwelvept
\else\fi

\theoremstyle{plain}
\newtheorem{Theorem}{Theorem}[section]

\newtheorem{Proposition}[Theorem]{Proposition}
\newtheorem{Lemma}[Theorem]{Lemma}
\newtheorem{Corollary}[Theorem]{Corollary}

\theoremstyle{definition}

\newtheorem{Definition}[Theorem]{Definition}
\newtheorem{Remark}[Theorem]{Remark}
\newtheorem{Construction}[Theorem]{Construction}
\newtheorem{Example}[Theorem]{Example}

\newtheorem{Notation}[Theorem]{Notation}

\newcommand{\ZZ}{\mathbf{Z}}

\newcommand{\RR}{\mathbb{R}}

\newcommand{\DD}{\mathbb{D}}

\newcommand{\NNNN}{\operatorname{N}}

\newcommand{\HH}{\operatorname{\mathcal{HH}}}

\newcommand{\EE}{\mathcal{E}}

\newcommand{\DDD}{\mathcal{D}}
\newcommand{\uni}{\mathbf{1}}
\newcommand{\CCC}{\mathcal{C}}

\newcommand{\PR}{\operatorname{Pr}^{\textup{L}}}
\newcommand{\PRT}{\operatorname{Pr}_{\textup{t}}^{\textup{L}}}
\newcommand{\PRTT}{\operatorname{Pr}_{\textup{t}+}^{\textup{L}}}
\newcommand{\PRTTT}{\operatorname{Pr}_{\textup{t}\pm}^{\textup{L}}}

\newcommand{\MMM}{\mathcal{M}}

\newcommand{\Hom}{\operatorname{Hom}}

\newcommand{\Ker}{\operatorname{Ker}}

\newcommand{\Spec}{\operatorname{Spec}}

\newcommand{\Perf}{\operatorname{Perf}}

\newcommand{\SP}{\operatorname{Sp}}

\newcommand{\Mod}{\operatorname{Mod}}

\newcommand{\SSS}{\mathcal{S}}
\newcommand{\wSSS}{\widehat{\mathcal{S}}}

\newcommand{\colim}{\operatorname{colim}}

\newcommand{\Cat}{\textup{Cat}_{\infty}}

\newcommand{\Map}{\operatorname{Map}}

\newcommand{\Fun}{\operatorname{Fun}}
\newcommand{\Alg}{\operatorname{Alg}}

\newcommand{\End}{\operatorname{End}}

\newcommand{\wCat}{\widehat{\textup{Cat}}_{\infty}}

\newcommand{\CAlg}{\operatorname{CAlg}}

\newcommand{\DEFOR}{\mathcal{D}ef}

\newcommand{\LM}{\operatorname{\mathcal{LM}}}

\newcommand{\hhh}{\operatorname{h}}
\newcommand{\Ind}{\operatorname{Ind}}

\newcommand{\assoc}{\operatorname{As}}

\newcommand{\eone}{\mathbf{E}_1}
\newcommand{\etwo}{\mathbf{E}_2}
\newcommand{\eenu}{\mathbf{E}_n}
\newcommand{\einf}{\mathbf{E}_\infty}

\newcommand{\LMod}{\operatorname{LMod}}
\newcommand{\RMod}{\operatorname{RMod}}

\newcommand{\HL}{\widehat{L}}

\newcommand{\Proof}{{\sl Proof.}\quad}
\newcommand{\QED}{{\unskip\nobreak\hfil\penalty50\quad\null\nobreak\hfil
{$\Box$}\parfillskip0pt\finalhyphendemerits0\par\medskip}}

\begin{document}

\title[Hochschild cohomology and Deformation]{Hochschild cohomology and Deformation theory of stable $\infty$-categories with t-structures}

\author{Isamu Iwanari}






\address{Mathematical Institute, Tohoku University, Sendai, Miyagi, 980-8578,
Japan}

\email{isamu.iwanari.a2@tohoku.ac.jp}

\maketitle

\section{Introduction}

Gerstenhaber \cite{Ger1} initiated the study of the deformation theory of associative algebras by means of Hochschild cohomology.
Hochschild cohomology carries a (graded) Lie algebra structure \cite{Ger2}.
This Lie algebra structure, together with its associative algebra structure, endows the Hochschild cohomology with the structure of a Gerstenhaber algebra.
This Gerstenhaber algebra structure governs the deformation theory of associative algebras.
Furthermore, the Hochschild cohomology of an algebra admits an enhancement from a Gerstenhaber algebra to an algebra over the operad of little $2$-disks,
i.e., an $\etwo$-algebra at the level of chain complexes.
This is commonly referred to as Deligne conjecture on Hochschild cohomology.
The conjecture has been solved through various approaches by McClure–Smith \cite{MS}, Kontsevich–Soibelman \cite{KS}, and Tamarkin \cite{Tam}.
Moreover, various generalized versions of the conjecture have been extensively studied, including those in categorical frameworks.

We turn to consider deformations of stable $\infty$-categories.
Let $k$ be a field.
Given a $k$-linear (presentable) stable $\infty$-category $\CCC$, one can define the Hochschild cochain complex
$\HH^\bullet(\CCC/k)$ which is defined as an $\etwo$-algebra.  
Let $\textup{Art}$ denote the $\infty$-category of Artin dg $k$-algebras.
Define $\DEFOR_{\CCC}:\textup{Art}\to \SSS$
be the deformation functor 
which sends $R$ to the space (i.e., $\infty$-groupoid) of deformations $(\DDD, \DDD\otimes_{R}k\simeq \CCC)$ of
$\CCC$, where $\DDD$ is an $R$-linear stable $\infty$-category, and $\SSS$ is the $\infty$-category of spaces.
The deformation theory of $\CCC$, or variants of compactly generated $\infty$-categories, has been investigated from several perspectives, for instance, in the setting of (derived) formal moduli problems \cite[X, Section 5]{DAG}, in the case of deformations to the ring of formal power series \cite{BKP}, and from the lifting problem of compact generators \cite{LV}.
In \cite{DAG}, the deformation theory of stable
$\infty$-categories is formulated in the local version of derived geometry over $\etwo$-algebras:
The deformation problem is defined as a functor $\DEFOR^{(2)}_{\CCC}:\textup{Art}_2\to \SSS$
from the $\infty$-category $\textup{Art}_2$ of Artin $\etwo$-algebras over $k$.
It was shown that, under restrictive conditions, $\DEFOR_{\CCC}^{(2)}$  is equivalent to the
$\etwo$-version of formal moduli problem $\mathcal{F}_{A\oplus \HH^\bullet(\CCC/k)}^{(2)}$ associated to the $\etwo$-algebra $\HH^\bullet(\CCC/k)$ while
there is always 
a map $\eta:\DEFOR_{\CCC}^{(2)}\to \mathcal{F}_{A\oplus \HH^\bullet(\CCC/k)}^{(2)}$ which is universal among maps to formal $\etwo$-moduli problems.
However, in general $\DEFOR_{\CCC}$ and $\DEFOR_{\CCC}^{(2)}$ do not satisfy the Schlessinger conditions of \cite[X]{DAG}
(see \cite{KL} for a specific example).
This implies that in general neither $\DEFOR_{\CCC}$ nor $\DEFOR_{\CCC}^{(2)}$ is a formal moduli problem.
Consequently, they can not be governed by the Hochschild cohomology or other algebraic invariants.
This has been a longstanding problem in the theory.
In other words, the deformation problem encoded by $\DEFOR_{\CCC}$ (or the likes) is somewhat naive,
and a more refined or appropriate formulation of the deformation problem is required.

The main purpose of this paper is to prove the following result (see Theorem~\ref{mainthm}, Theorem~\ref{mainthm2} for precise statements):

\begin{Theorem}
\label{intromain1}
Let $\EE$ be a $k$-linear presentable stable $\infty$-category equipped with a left and right complete, accessible $t$-structure $(\EE^{\le0},\EE^{\ge0})$.
Let $\textup{Deform}_{\EE}:\textup{Art}_2\to \wSSS$ be the deformation functor which carries $R\in \textup{Art}_2$ to the space $\textup{Deform}_{\EE}(R)$
consisting of deformations 
\[
\big{(}(\DDD,\DDD^{\le0}), (\DDD\otimes_{R}k,(\DDD\otimes_{R}k)^{\le0})\simeq (\EE, \EE^{\le0})\big{)}
\]
of $(\EE,\EE^{\le0})$ where  $(\DDD,\DDD^{\le0})$ indicates an $R$-linear presentable $\infty$-category equipped with a left and right complete $t$-structure (see Definition~\ref{deffunctorcomplete} for the detail). 
Then the deformation functor $\textup{Deform}_{\EE}$
is a formal $\etwo$-moduli problem $\mathcal{F}_{k\oplus \HH^\bullet(\EE/k)}^{(2)}$ (in Lurie's sense) that corresponds to the augmented $\etwo$-algebra $k\oplus \HH^\bullet(\EE/k)\to k$. 
In particular, there is a canonical equivalence of spaces
\[
\textup{Deform}_{\EE}(R)\simeq \Map_{\Alg_2(\Mod_k)}(\DD_2(R), \HH^\bullet(\EE/k))
\]
for $R\in \textup{Art}_2$. Here $\Alg_2(\Mod_k)$ is the $\infty$-category of $\etwo$-algebras over $k$, and $\DD_2$ is the $\etwo$-Koszul duality functor.
\end{Theorem}

If we focus on deformations to commutative bases, we have:

\begin{Corollary}
Let $\textup{Art}$ be the $\infty$-category of Artin dg $k$-algebras. 
Let $\textup{Deform}_{\EE}^\infty:\textup{Art}\to \textup{Art}_2\to \wSSS$ be
the composite of the forgetful functor $\textup{Art}\to \textup{Art}_2$ and $\textup{Deform}_{\EE}: \textup{Art}_2\to \wSSS$.
Then the deformation functor $\textup{Deform}^\infty_{\EE}$
is a formal $\mathbf{E}_{\infty}$-moduli problem.
When $k$ is a field of characteristic zero,
$\textup{Deform}^\infty_{\EE}$ is equivalent to
the formal moduli problem $\mathcal{F}^{(\infty)}_{\HH^\bullet(\EE/k)[1]}$ corresponding to the underlying dg Lie algebra $\HH^\bullet(\EE/k)[1]$
associated to 
the $\etwo$-algebra $\HH^\bullet(\EE/k)$.
Namely, $\mathcal{F}^{(\infty)}_{\HH^\bullet(\EE/k)[1]}$ is defined by $\mathcal{F}^{(\infty)}_{\HH^\bullet(\EE/k)[1]}(R)=\Map_{Lie_k}(\DD_\infty(R), \HH^\bullet(\EE/k)[1])$
where $\DD_\infty:((\CAlg_k)_{/k})^{op}\to Lie_k$ is the Koszul duality functor which is a right adjoint to
the Chevalley-Eilenberg cochain functor $Ch^\bullet:Lie_k \to ((\CAlg_k)_{/k})^{op}$. 
\end{Corollary}

As a direct application, we deduce an obstruction theory for the deformation theory (see Theorem~\ref{mainthm4}).

\vspace{2mm}

We note that left complete 
$t$-structures (cf. \cite{HA}) represent a natural form of completeness in practice.
Many important classes of stable $\infty$-categories naturally admit left complete 
$t$-structures—for example, the stable 
$\infty$-categories of quasi-coherent sheaves/complexes on derived algebraic stacks, of 
$D$-modules, of connective differential graded algebras, and so on.
Moreover, the process of left completion sends any $t$-structure to a left complete one.

We explain Theorem~\ref{intromain1} from the perspective of Koszul duality.
Our main result may be understood as
a realization of
categorified Koszul duality.
By this, we mean 
$\eone$-Koszul duality at the level of $\infty$-categories, where
a module $M$ over an associative algebra $A$ is replaced 
by an $\mathcal{A}$-linear stable $\infty$-category $\MMM$
with a monoidal stable $\infty$-category $\mathcal{A}$.
Indeed, 
Lurie's construction 
of $\eta:\DEFOR_{\CCC}^{(2)}\to \mathcal{F}_{A\oplus \HH^\bullet(\CCC/k)}^{(2)}$ can be regarded as Koszul dual construction at level of stable $\infty$-categories, although this aspect is not explicitly emphasized.
In this paper, we build on this viewpoint by considering an ``adjoint" version of the Koszul dual construction, taking into account left complete $t$-structures.
(A reflection of categorified Koszul duality can also be found in \cite[Chapter 8]{AL}.)
In classical Koszul duality, a completeness condition is essential for establishing an equivalence between modules over an algebra and modules over its Koszul dual.
Analogously, in the categorified setting, issues of completeness naturally arise, and what we found is that
the notion of left complete $t$-structures
plays a key role.

There are several possible directions in which the results and methodology of this study could be extended.
For instance, let us consider a $k$-linear presentable $\infty$-category $\CCC^\otimes$
endowed with an $\eenu$-monoidal structure and a left complete $t$-structure compatible with the monoidal structure.
The deformation theory of $\CCC^\otimes$ on the level of the $\eenu$-monoidal $\infty$-category with the left complete $t$-structure 
seems likely to be governed by the $\mathbf{E}_{n+1}$-Hochschild cochain complex (see \cite{FR}) in the setting of formal $\mathbf{E}_{n+2}$-moduli problems.

\vspace{2mm}

This paper is organized as follows.
In Section~\ref{Pre} we fix some notation and terminology
and review notions and some results.
In Section~\ref{adjointsec} we discuss left completions, which appear in categorified Koszul duality.
Along the way we also obtain an application to deformations of stable presentable $\infty$-categories without $t$-structures (Corollary~\ref{without}).
In Section~\ref{mainsec}
we prove the main results of this paper.

\vspace{2mm}

{\it Acknowledgements.}
This work was partially supported by JSPS KAKENHI grant. 
The author would like to thank T. Matsuoka for the discussion about Koszul
duality and earlier versions of this paper.

\section{Preliminaries}
\label{Pre}

{\it Notation and Terminology.}
We use 
the language of $(\infty,1)$-categories.
We use the theory of {\it quasi-categories} as a model of $(\infty,1)$-categories.
Our main references are \cite{HTT}
 and \cite{HA} as well as subsequent works \cite{Kero} and \cite{RV}.
Following \cite{HTT}, we shall refer to quasi-categories
as {\it $\infty$-categories}.
To an ordinary category, we can assign an $\infty$-category by taking
its nerve, and therefore
when we treat ordinary categories we often omit the nerve $\NNNN(-)$
and directly regard them as $\infty$-categories.
Here is a list of some of the notation and conventions we frequently use:

\begin{itemize}

\item $\simeq$: categorical equivalence or homotopy equivalence.  The symbol $\simeq$ is used to indicate isomorphisms and equivalences in categories and
$\infty$-categories.

\item $\CCC^{op}$: the opposite $\infty$-category of an $\infty$-category.

\item $\Map_{\mathcal{C}}(C,C')$: the mapping space from an object $C\in\mathcal{C}$ to $C'\in \mathcal{C}$ where $\mathcal{C}$ is an $\infty$-category.

\item $\Fun(A,B)$: the function complex for simplicial sets $A$ and $B$. If $A$ and $B$ are $\infty$-categories, we regard $\Fun(A,B)$ as the functor category. 

\item $\SSS$: the $\infty$-category of small spaces ($\infty$-groupoids). We write $\wSSS$ for the $\infty$-category of spaces in an enlarged universe. 

\item $\Cat$: the $\infty$-category of small $\infty$-categories. We write $\wCat$ for the $\infty$-category of $\infty$-categories in an enlarged universe. 

\item $\Ind(\CCC)$: the Ind-category of $\CCC$.

\end{itemize}

We will use the theory of higher algebras.
We refer the reader to \cite{HA}.
Here is a list of (some) of the notation about $\infty$-operads and algebras over them that we will use:

\begin{itemize}
\label{notation2}

\item $\mathbf{E}^\otimes_n$: the $\infty$-operad of little $n$-cubes (cf. \cite[5.1]{HA}).

\item $\Alg_n(\MMM)$:
For a symmetric monoidal $\infty$-category $\MMM^\otimes$,
we write $\Alg_{n}(\MMM)$ for
the $\infty$-category of algebra objects over $\mathbf{E}^\otimes_n$ in $\MMM^\otimes$.
We refer to an object of $\Alg_{n}(\MMM)$ as an $\eenu$-algebra in $\MMM$.
If we denote by $\assoc^\otimes$ the associative operad (\cite[4.1.1]{HA}),
there is the standard equivalence $\assoc^\otimes\simeq \eone^\otimes$
of $\infty$-operads. We usually identify $\Alg_{1}(\MMM)$
with the $\infty$-category $\Alg_{\assoc}(\MMM)$, that is, the 
$\infty$-category of associative algebras in $\MMM$.

\item $\CAlg(\mathcal{M})$: the $\infty$-category of commutative
algebra objects in a symmetric
monoidal $\infty$-category $\mathcal{M}^\otimes$.
Instead of $\Alg_{\einf}(\MMM)$ or $\Alg_{\textup{comm}}(\MMM^\otimes)$
we write $\CAlg(\mathcal{M})$ for the $\infty$-category of algebra objects over the $\einf$-operad or the commutative operad in $\MMM^\otimes$.

\item $\CAlg_k$: $\CAlg_k=\CAlg(\Mod_k)$.

\item $\Mod_R(\mathcal{M})$: 
the $\infty$-category of
$R$-module objects,
where $\mathcal{M}^\otimes$
is a symmetric monoidal $\infty$-category and $R\in \CAlg(\MMM)$.

\item $\Mod_R$: When $\MMM^\otimes$ is the symmetric monoidal $\infty$-category $\SP^\otimes$
of spectra, 
we write $\Mod_R^\otimes$
for
$\Mod_R^\otimes(\SP)$.
We denote by $\Mod_R$ the underlying category, that is, $\Mod_R(\SP)$.

\item 
$\LMod(\MMM)$:
the $\infty$-category $\Alg_{\LM}(\MMM)$ of $\LM^\otimes$-algebras
over $\LM^\otimes$: the $\infty$-operad defined in \cite[4.2.1.7]{HA}.
We may consider an object of  $\LMod(\MMM)$ to be a pair 
of $(A,M)$ such that $A$ is an associative algebra object in $\MMM$
and a left $A$-module $M$.
There is a Cartesian fibration $\LMod(\MMM)\to \Alg_{\assoc}(\MMM)\simeq \Alg_1(\MMM)$
which sends $(A,M)$ to $A$.
Another canonical morphism is the forgetful functor $\LMod(\MMM)\to \MMM$
which is informally given by $(A,M)\mapsto M$.

\item $\LMod_A(\MMM)$: the $\infty$-category of left $A$-module objects in $\MMM$.
For $A\in \Alg_1(\MMM)$, we define $\LMod_A(\CCC)$
to be the fiber of $\LMod(\CCC)\to \Alg_1(\CCC)$ over $A$ in $\Cat$.

\item $\RMod(\MMM)$: the right module version of $\LMod(\MMM)$

\item $\RMod_B(\MMM)$: the $\infty$-category of right $B$-module objects in $\MMM$.

\end{itemize}

Let $k$ denote a fixed field, which serves as the base field.
Unless otherwise stated, we adopt the cohomological index.

\subsection{Generalities on $t$-structures}

Let us recall the notion of $t$-structures on stable $\infty$-categories.
Let $\CCC$ be a stable $\infty$-category.
Recall that the homotopy category $\hhh(\CCC)$ has the structure of a triangulated category.
A $t$-structure on $\CCC$ is a pair of full subcategories $(\CCC^{\le0},\CCC^{\ge0})$ of $\CCC$
such that $(\hhh(\CCC^{\le0}), \hhh(\CCC^{\ge0}))$ is a $t$-structure on the triangulated category $\hhh(\CCC)$. The pair
is required to satisfy the following conditions:
If we write $\CCC^{\le n}$ and $\CCC^{\ge n}$ for $\CCC^{\le 0}[-n]$ and $\CCC^{\ge 0}[-n]$, respectively, then
\begin{enumerate}
\item $\CCC^{\le-1}\subset \CCC^{\le0}$ and $\CCC^{\ge1}\subset \CCC^{\ge0}$.

\item $\Hom_{\hhh(\CCC)}(C,C')=0$ for $C\in \CCC^{\le0}$ and $C'\in \CCC^{\ge1}$.

\item For any $C\in \CCC$, there exists a cofiber square 
\[
\xymatrix{
D \ar[r] \ar[d] & C \ar[d] \\
0 \ar[r] &  E
}
\]
such that $D\in \CCC^{\le0}$ and $E\in \CCC^{\ge1}$.
\end{enumerate}

Suppose that $(\CCC^{\le0},\CCC^{\ge0})$ is a $t$-structure.
The inclusion $\CCC^{\le n}\hookrightarrow\CCC$ has a right adjoint 
\[
\tau^{\le n}:\CCC\longrightarrow \CCC^{\le n},
\]
and the inclusion $\CCC^{\ge n}\hookrightarrow\CCC$
has a left adjoint 
\[
\tau^{\ge n}:\CCC\longrightarrow \CCC^{\ge n}.
\]

There exists a canonical equivalence $\tau^{\ge 0}\circ \tau^{\le0}\simeq \tau^{\le 0}\circ \tau^{\ge0}$.
We define
\[
H^0:\CCC\longrightarrow \CCC^{\le0}\cap \CCC^{\ge0}
\]
to be $\tau^{\ge 0}\circ \tau^{\le0}$ (equivalently, $\tau^{\le 0}\circ \tau^{\ge0}$).
The intersection $\CCC^{\le0}\cap \CCC^{\ge0}$ is called the heart.
The heart $\CCC^{\le0}\cap \CCC^{\ge0}$ is an abelian category (see \cite{BBD}).
We write $\CCC^{\heartsuit}=\CCC^{\le0}\cap \CCC^{\ge0}$.
We define $H^n=H^0\circ [n]:\CCC\to \CCC^{\le0}\cap \CCC^{\ge0}$ where 
$[n]$ is the $n$-fold shift (suspension) functor.

\begin{Remark}
Given a cofiber sequence $C\to C'\to C''$ in $\CCC$, there exists a canonical long exact sequence
\[
\cdots \to H^n(C)\to H^{n}(C')\to H^n(C'')\to H^{n+1}(C)\to \cdots
\]
in $\CCC^\heartsuit$ (see \cite{BBD}).
\end{Remark}

We often use the notation $(\CCC^{\le0}, \CCC^{\ge0})$ to denote a $t$-structure on $\CCC$.
Note that $\CCC^{\ge0}$ is uniquely determined by $\CCC^{\le0}$: the full subcategory $\CCC^{\ge0}$ is spanned by those objects $C'$ such that $\Map_{\CCC}(C,C')$ is contractible for any $C\in \CCC^{\le -1}$. Dually, $\CCC^{\le0}$ is uniquely determined by $\CCC^{\ge0}$.
By $(\CCC,\CCC^{\le0})$ we mean a stable $\infty$-category endowed with the $t$-structure determined by a full subcategory $\CCC^{\le0}$.

\begin{Definition}
Let $(\CCC^{\le 0}, \CCC^{\ge 0})$ be a $t$-structure on a stable $\infty$-category $\CCC$.
\begin{enumerate}
\item For any $n \in \ZZ$, the left adjoint functor $\tau^{\ge n}: \CCC \to \CCC^{\ge n}$ to the inclusion
$\CCC^{\ge n} \subset \CCC$ induces a morphism $C \to \tau^{\ge n}(C)$, determined by the unit map. For all $n \in \ZZ$, this yields the sequence
\[
C \to \cdots \to \tau^{\ge n-1}(C) \to \tau^{\ge n}(C) \to \tau^{\ge n+1}(C) \to \cdots
\]
in $\CCC$.
The $t$-structure is left complete if, for any $C \in \CCC$, the induced morphism $C \to \varprojlim_{n \in \ZZ} \tau^{\ge n}(C)$ is an equivalence.
\item The right adjoint functor $\tau^{\le n}: \CCC \to \CCC^{\le n}$ to the inclusion
$\CCC^{\le n} \subset \CCC$ induces a morphism $\tau^{\le n}(C) \to C$, determined by the counit map. For all $n \in \ZZ$, this yields the sequence
\[
\cdots \to \tau^{\le n-1}(C) \to \tau^{\le n}(C) \to \tau^{\le n+1}(C) \to \cdots \to C
\]
in $\CCC$.
The $t$-structure is right complete if, for any $C \in \CCC$, the induced morphism $\varinjlim_{n \in \ZZ} \tau^{\le n}(C) \to C$ is an equivalence.
\item The $t$-structure is left separated if $\bigcap_{n \in \ZZ} \CCC^{\le n} = 0$.
Dually, it is right separated if $\bigcap_{n \in \ZZ} \CCC^{\ge n} = 0$.
\item If $\CCC$ admits filtered colimits and the full subcategory $\CCC^{\ge 0}$ is closed under filtered colimits in $\CCC$, we say that the $t$-structure is compatible with filtered colimits.
\item The $t$-structure $(\CCC^{\le 0}, \CCC^{\ge 0})$ is accessible if the full subcategory $\CCC^{\le 0}$ is presentable.
\item Suppose that $\CCC$ has a monoidal structure. If the tensor product functor $\CCC \times \CCC \to \CCC$ sends $\CCC^{\le 0} \times \CCC^{\le 0}$ to $\CCC^{\le 0}$ and the unit object lies in $\CCC^{\le 0}$, we say that the $t$-structure is compatible with the monoidal structure.
\end{enumerate}
\end{Definition}

\begin{Remark}
If $C\in \cap_{n\in \ZZ}\CCC^{\le n}$, then $\tau^{\ge n}C\simeq0$ for all $n\in \ZZ$.
It follows that if a $t$-structure $(\CCC^{\le0}, \CCC^{\ge0})$ is left complete,
then it is left separated.
Dually, if a $t$-structure is right complete,
then it is right separated.
\end{Remark}

\begin{Example}
\label{connectivetstr}
Let $A$ be a connective associative ring spectrum, and let $\LMod_A$ be the stable $\infty$-category of left $A$-module spectra.
Define $\LMod_A^{\ge 0}$ as the full subcategory of $\LMod_A$ spanned by objects $M$ such that $\pi_n(M) = 0$ for all $n < 0$, and $\LMod_A^{\le 0}$ as the full subcategory spanned by objects $M$ such that $\pi_n(M) = 0$ for all $n > 0$.
Then $(\LMod_A^{\le 0}, \LMod_A^{\ge 0})$ defines an accessible, left complete, and right complete $t$-structure (see \cite[7.1.1.13]{HA}).
The full subcategory $\LMod_A^{\le 0}$ is also characterized as the smallest full subcategory which contains $A$ and is closed under small colimits.
Moreover, the $t$-structure is compatible with filtered colimits, as $\LMod_A^{\ge 0} \subset \LMod_A$ is closed under filtered colimits.
This $t$-structure is compatible with the monoidal structure when $A$ is an $\etwo$-algebra.
The heart $\LMod_A^{\heartsuit}$ is equivalent to the abelian category of left $\pi_0(A)$-modules.
\end{Example}

\begin{Definition}
Let $A$ be an $\eenu$-algebra in $\Mod_k$.
By convention, an $\mathbf{E}_{0}$-algebra means an object of $\Mod_k$.
We say that $A$ is $m$-coconnective if
$H^i(A)=0$ for any $i<0$ and $0<i<m$, and the canonical map $k\to H^0(A)$ is an isomorphism.
When $A$ is $1$-coconnective, we say simply that $A$ is coconnective.
We say that $A$ is locally finite if $H^i(A)$ is finite dimensional for any $i\in \ZZ$.
\end{Definition}

\begin{Proposition}
\label{tstrcoconnective}
Let $D$ be an $\eone$-algebra in $\Mod_k$.
Assume that $D$ is coconnective.
Let $\LMod_D^{\le 0}$ be the smallest full subcategory of $\LMod_D$ containing $D$ and closed under small colimits and extensions.
Let $\LMod_D^{\ge 0}$ be the full subcategory spanned by objects $M$ such that their image under the forgetful functor to $\Mod_k$ lies in $\Mod_k^{\ge 0}$, i.e., $H^i(M) = 0$ for all $i < 0$.
Then $(\LMod_D^{\le 0}, \LMod_D^{\ge 0})$ defines an accessible and right complete $t$-structure.
Moreover, it is compatible with filtered colimits.
\end{Proposition}

\Proof
When $D$ is commutative, the assertion is proved in \cite[VIII, 4.5.4]{DAG}.
The proof also works for associative algebras.
\QED

\begin{Remark}
The $t$-structure $(\LMod_D^{\le0}, \LMod_D^{\ge0})$
is not left complete in general (cf. \cite[VIII, 4.5.5]{DAG})
\end{Remark}

\begin{Remark}
When $D$ is (promoted to) an $\etwo$-algebra, then $\LMod_D$ has the canonical
monoidal structure. From the definition of 
$\LMod_D^{\le0}$ we see that $\LMod_D^{\le0}$ is closed under tensor products.
\end{Remark}

\begin{Remark}
\label{unifiedtstr}
Example~\ref{connectivetstr} and Proposition~\ref{tstrcoconnective} are unified as follows.
For $A\in \Alg_1(\SP)$ we define $\LMod_A^{\le0}$ to be the smallest full subcategory
which contains $A$ and is closed under small colimits and extensions. By \cite[1.4.4.11]{HA},
$\LMod_A^{\le0}$ determines an accessible $t$-structure on $\LMod_A$ such that $\LMod_A^{\ge0}$
consists of those objects $M$ such that $\pi_0(\Map_{\LMod_A}(A[n],M))=\pi_n(M)=H^{-n}(M)=0$ for $n>0$.
The full subcategory $\LMod_A^{\ge0}$ is closed under filtered colimits, and the $t$-structure is right complete.
When $A$ is an $\etwo$-algebra, the $t$-structure is compatible with the monoidal structure on $\LMod_A$.
As seen in Example~\ref{connectivetstr} and Proposition~\ref{tstrcoconnective}, the property of the resulting $t$-structure
depends on $A$.
\end{Remark}

\begin{Example}
Let $\CCC$ be a small stable $\infty$-category.
Let $(\CCC^{\le0}, \CCC^{\ge0})$ be a $t$-structure on $\CCC$.
Then $(\Ind(\CCC^{\le0}),\Ind(\CCC^{\ge0}))$ defines a $t$-structure on $\Ind(\CCC)$.
By the construction, the $t$-structure is compatible with filtered colimits.
If the $t$-structure on $\CCC$ is right bounded, then 
$(\Ind(\CCC^{\le0}),\Ind(\CCC^{\ge0}))$ is right complete (see \cite[VIII, Lemma 5.4.1]{DAG}).
\end{Example}

\begin{Example}
Let $X$ be a derived/spectral scheme (or more generally, derived/spectral Deligne-Mumford stack that satisfies a mild condition).
Then the stable $\infty$-category $\textup{QCoh}(X)$ of quasi-coherent sheaves on $X$ has a $t$-structure
which is accessible, left complete, and right complete. It is also compatible with filtered colimits and the monoidal structure (see \cite[VIII, 2.3.18]{DAG}).
\end{Example}

\begin{Definition}
Let $(\CCC,\CCC^{\le0})$ and $(\DDD,\DDD^{\le0})$ be stable $\infty$-categories endowed with $t$-structure.
\begin{enumerate}
\item An exact functor $F:\CCC\to \DDD$ is said to be right $t$-exact if the essential image of $\CCC^{\le0}$ is contained in $\DDD^{\le0}$.

\item An exact functor $F:\CCC\to \DDD$ is said to be left $t$-exact if the essential image of $\CCC^{\ge0}$ is contained in $\DDD^{\ge0}$.

\item An exact functor $F:\CCC\to \DDD$ is said to be $t$-exact if $F$ is right $t$-exact and left $t$-exact.

\end{enumerate}
\end{Definition}

Let $\PR_{\mathbb{S}}$ be the $\infty$-category of presentable stable $\infty$-categories (cf. \cite[Definition 5.5.3.1]{HTT}, \cite{HA}). A morphism $\CCC\to \DDD$ in $\PR_{\mathbb{S}}$ is a functor that preserves small colimits.
Let $\PRT$ be the $\infty$-category of presentable stable $\infty$-categories equipped with accessible $t$-structures.
An object of $\PRT$ is a pair $(\CCC,\CCC^{\le0})$ where $\CCC$ is a presentable stable $\infty$-category, and $\CCC^{\le0}$ is a  full subcategory
inducing a $t$-structure,
such that $\CCC^{\le0}$ itself is presentable.
A morphism $(\CCC, \CCC^{\le0})\to (\DDD,\DDD^{\le0})$ is a right $t$-exact functor which preserves small colimits.

Let $\PRTT$ be the full subcategory of $\PRT$, spanned by presentable stable $\infty$-categories equipped with accessible right complete $t$-structures.
Let $\PRTTT$ be the full subcategory of $\PRT$, spanned by presentable stable $\infty$-categories equipped with accessible left complete and right complete $t$-structures.
There exists an adjoint pair
\[
\xymatrix{
I^+: \PRTT \ar@<0.5ex>[r] &  \PRT:\widehat{R}  \ar@<0.5ex>[l]  
}
\]
such that $I^+$ is the inclusion, and the right adjoint $\widehat{R}$ carries $(\CCC,\CCC^{\le0})$
to the right completion 
\[
\widehat{R}(\CCC,\CCC^{\le0})=(\colim_{n\to \infty} \CCC^{\le n},\CCC^{\le0}) 
\]
where $\colim_{n\to \infty} \CCC^{\le n}$ is the colimit of the sequence of inclusions
$\CCC^{\le 0}\to \CCC^{\le 1}\to \CCC^{\le 2}\to \cdots\to \CCC^{\le n}\to \cdots$ in $\PR_{\mathbb{S}}$
(cf. \cite[5.5.3.4, 5.5.3.18]{HTT}). 

The left completion $\HL((\CCC,\CCC^{\le0}))$ of $(\CCC,\CCC^{\le0})$
is defined as
$\lim_{n\to -\infty}\CCC^{\ge n}$ in $\PR_{\mathbb{S}}$
of the sequence 
\[
\cdots \to \CCC^{\ge n-1}\to \CCC^{\ge n}\to \CCC^{\ge n+1} \to \cdots
\]
induced by truncations.
We usually write $\HL(\CCC)$ for $\HL((\CCC,\CCC^{\le0}))$.
If $(\CCC,\CCC^{\le0})\in \PRTT$, then the left completion $\HL(\CCC)$ is also accessible and right complete (see \cite[VIII, 4.6.13]{DAG}).
There exists an adjoint pair
\[
\xymatrix{
\HL :\PRTT \ar@<0.5ex>[r] &  \PRTTT:I^{\pm}  \ar@<0.5ex>[l]  
}
\]
such that $I^\pm$ is the inclusion, and $\HL$ is the left adjoint given by the left completion.

{\it Symmetric monoidal structures.}
We recall the symmetric monoidal structure on $\PR_{\mathbb{S}}$ (see \cite[4.8]{HA}).
Let $\mathcal{B}$ and $\CCC$ be presentable $\infty$-categories. 
The tensor product of $\mathcal{B}$ and $\CCC$
is defined as an object
$\mathcal{B}\otimes \CCC$ in $\PR_{\mathbb{S}}$ endowed with a functor  
$\mathcal{B}\times \CCC\to \mathcal{B}\otimes \CCC$
which preserves small colimits separately in each variable such that for any $\DDD\in \PR_{\mathbb{S}}$ the composition gives the fully faithful functor
\[
\Fun^\textup{L}(\mathcal{B}\otimes \CCC,\DDD)\to \Fun(\mathcal{B}\times \CCC,\DDD)
\]
whose essential image $\Fun'(\mathcal{B}\times \CCC,\DDD)$ consists of functors 
which preserves small colimits separately in each variable.
Indeed, $\mathcal{B}\otimes \CCC$ is equivalent to the full subcategory $\Fun^{\textup{R}}(\mathcal{B}^{op}, \CCC)$
of $\Fun(\mathcal{B}^{op}, \CCC)$ spanned by limit-preserving functors. This $\infty$-category
$\Fun^{\textup{R}}(\mathcal{B}^{op}, \CCC)$ is presentable (\cite[4.8.1.16]{HA}).

Next, we recall the symmetric monoidal structure on $\PRT$ (see \cite[VIII, 4.6.1]{DAG}).
Let $(\mathcal{B},\mathcal{B}^{\le0})$ and $(\CCC,\CCC^{\le0})$ be presentable $\infty$-categories endowed with accessible $t$-structures. 
The tensor product of $(\mathcal{B},\mathcal{B}^{\le0})$ and $(\CCC,\CCC^{\le0})$
is defined as
$\mathcal{B}\otimes \CCC$ (in $\PR_{\mathbb{S}}$) 
together with the smallest full subcategory $(\mathcal{B}\otimes \CCC)^{\le0}$
which contains the essential image of $\mathcal{B}^{\le0}\times \CCC^{\le0}\to \mathcal{B}\otimes \CCC$
and is closed under small colimits and extensions.
By \cite[1.4.5.11]{HA}, $(\mathcal{B}\otimes \CCC)^{\le0}$ defines an accessible $t$-structure.
The universal property of $(\mathcal{B}\otimes \CCC, (\mathcal{B}\otimes \CCC)^{\le0})$ is described as follows: For any $(\DDD, \DDD^{\le0})\in \PRT$
the composition with the canonical functor $\mathcal{B}\times \CCC\to \mathcal{B}\otimes \CCC$ gives an equivalence
\[
\Fun_{\PRT}((\mathcal{B}\otimes \CCC, (\mathcal{B}\otimes \CCC)^{\le0}),(\DDD,\DDD^{\le0}))\to \Fun^{\le}(\mathcal{B}\times \CCC,\DDD)
\]
where $\Fun^{\le}(\mathcal{B}\times \CCC,\DDD)$ is the full subcategory of $\Fun'(\mathcal{B}\times \CCC,\DDD)$,
which consists of those functors $F$ such that the essential image of $\mathcal{B}^{\le0}\times \CCC^{\le0}$ 
is contained in $\DDD^{\le0}$.

If $(\mathcal{B},\mathcal{B}^{\le0})$ and $(\CCC,\CCC^{\le0})$ are right complete,
then the tensor product $(\mathcal{B}\otimes \CCC, (\mathcal{B}\otimes \CCC)^{\le0})$ in $\PRT$
is also right complete (see \cite[VIII, 4.6.11]{DAG}).
Thus, $\PRTT$ inherits a symmetric monoidal structure from that of $\PRT$.
The left adjoint $I^+$ is promoted
to a symmetric monoidal functor.

\begin{Example}
Let $A$ and $B$ be associative ring spectra, that is, $\eone$-algebras in $\SP$.
The stable $\infty$-catgories $\LMod_A$ and $\LMod_B$ are equipped with $t$-structures described in Remark~\ref{unifiedtstr}.
Then $\LMod_A\otimes \LMod_B$ has the $t$-structure determined by the full subcategory
$(\LMod_A\otimes \LMod_B)^{\le0}$ corresponding to the smallest full subcategory of $\LMod_{A\otimes B}$ which contains $A\otimes B$
and is closed under small colimits and extensions through the canonical equivalence $\LMod_A\otimes \LMod_B\simeq \LMod_{A\otimes B}$ (cf. Section~\ref{MOA}).

Let $C\to k$ be an augmented $\etwo$-algebra in $\Mod_k$.
The stable $\infty$-catgories $\Mod_k$ and $\LMod_C$ are equipped with $t$-structures described in Remark~\ref{unifiedtstr}.
Then the relative tensor product $\Mod_k\otimes_{\LMod_A}\Mod_k$ has the $t$-structure determined by
the smallest full subcategory $(\Mod_k\otimes_{\LMod_A}\Mod_k)^{\le0}$, which contains the essential image of the functor 
$\Mod_k^{\le0}\times\Mod_k^{\le0}\to \Mod_k\otimes_{\LMod_A}\Mod_k$ (arising from the bar construction)
and is closed under small colimits and extensions.
Observe that through the canonical equivalence $\Mod_k\otimes_{\LMod_A}\Mod_k\simeq \LMod_{k\otimes_A k}$ (cf. Lemma~\ref{moduleproduct}),
$(\Mod_k\otimes_{\LMod_A}\Mod_k)^{\le0}$ corresponds to the smallest full subcategory of $\LMod_{k\otimes_A k}$,
which contains $k\otimes_Ak$ and is closed under small colimits and extensions.
\end{Example}

The $\infty$-category $\PRTTT$ also has a symmetric monoidal structure.
Let $(\mathcal{B},\mathcal{B}^{\le0})$ and $(\CCC,\CCC^{\le0})$ be objects of $\PRTTT$.
Explicitly, the tensor product of $(\mathcal{B},\mathcal{B}^{\le0})$ and $(\CCC,\CCC^{\le0})$ in $\PRTTT$
is given by 
\[
\HL((\mathcal{B}\otimes\CCC,  (\mathcal{B}\otimes \CCC)^{\le0})).
\]
We usually write $\HL(\mathcal{B}\otimes\CCC)$ or $\mathcal{B}\widehat{\otimes}\CCC$ for it.
The functor $\HL:\PRTT\to \PRTTT$ is promoted to a symmetric monoidal functor. 
Since $\HL$ is a left adjoint functor, $\HL$ commutes with the relative tensor products.
The inclusion $I^{\pm}$
is not a symmetric monoidal functor but a lax symmetric monoidal functor.

\subsection{Modules over algebras.}
\label{MOA}
Let $A$ be a commutative ring spectrum.
We write $\PR_A=\Mod_{\Mod_A}(\PR_{\mathbb{S}})$.
For $B\in \Alg_{1}(\Mod_A)$, we denote by $\LMod_B(\Mod_A)$ (resp. $\RMod_A(\Mod_A)$)
(or simply $\LMod_B$ and (resp. $\RMod_B$)) the $\infty$-category of left $B$-modules
(resp. right $B$-module spectra) (there is a canonical equivalence $\LMod_B(\Mod_A)\simeq \LMod_B(\SP)$ induced by the forgetful functor $\Mod_A\to \SP$).
As constructed \cite[4.8.5]{HA}, there is a symmetric monoidal functor
\[
\Theta_A:\Alg_{1}(\Mod_A)\longrightarrow (\PR_A)_{\Mod_A/}
\]
which carries $B$ to $\alpha:\Mod_A\to \LMod_{B}(\Mod_A)$
such that $\alpha$ is determined by $\SP\stackrel{\otimes A}{\to} \Mod_A\to \LMod_{B}(\Mod_A)$ which carries the sphere spectrum
to $B$. The homomorphism $f:B\to C$ maps to the base change functor $f^*:\LMod_B\to \LMod_C$.
Here the subscript in $(\PR_A)_{\Mod_A/}$ means the undercategory.
The functor $\Theta_A$ is a fully faithful left adjoint functor, see \cite[4.8.5.11]{HA}.
The right adjoint $E_A:(\PR_A)_{\Mod_A/}\to \Alg_{1}(\Mod_A)$
sends $\alpha:\Mod_A\to \CCC$ to the opposite algebra of the endomorphism algebra of $\alpha(A)\in \CCC$, defined as an algebra object of $\Mod_A$. 
Let $\iota_A:\Alg_1(\Mod_A)\to \Alg_1(\Mod_A)$
be the equivalence which carries $B$ to the opposite algebra $B^{op}$.
This equivalence is induced by the automorphism of the $\eone$-operad $\eone^\otimes\to \eone^\otimes$
given by the reflection of lines ($x\mapsto -x$).
We have
$\Theta_A\circ \iota_A:\Alg_{1}(\Mod_A)\longrightarrow (\PR_A)_{\Mod_A/}$
which carries $B$ to $\alpha:\Mod_A\to \RMod_{B}(\Mod_A)=\LMod_{B^{op}}(\Mod_A)$.

Since $\Theta_A$ is a symmetric monoidal functor, 
by Dunn additivity theorem \cite[5.1.2]{HA}, it determines the fully faithful symmetric monoidal left adjoint functor 
\[
\Alg_1(\Theta_A):\Alg_{2}(\Mod_A)\simeq \Alg_{1}(\Alg_{1}(\Mod_A))\to  \Alg_{1}((\PR_A)_{\Mod_A/})\simeq \Alg_{1}(\PR_A).
\]
The right adjoint of $\Alg_1(\Theta_A)$ is $\Alg_1(E_A)$ induced by the lax symmetric monoidal functor $E_A$.
The induced right adjoint functor $\Alg_1(E_A):\Alg_{1}(\PR_A)\to \Alg_{2}(\Mod_A)$ carries $\MMM^\otimes\in \Alg_{1}(\PR_A)$ to (the opposite algebra of) the endomorphism algebra $\End(\mathsf{1}_{\MMM})$
of the unit object $\mathsf{1}_{\MMM}$, defined as an $\etwo$-algebra object in $\Mod_A$.
By abuse of notation, we often write $\Theta_A$ and $E_A$ for $\Alg_1(\Theta_A)$ and $\Alg_1(E_A)$, respectively.

For $B\in \Alg_2(\Mod_A)$, since the underlying $A$-linear presentable $\infty$-category
of $\Theta_A(B)$ is $\LMod_B$, $\LMod_B$ is promoted to an object of $\Alg_1(\PR_A)$.
In particular, $\LMod_B$ is equipped with the structure of an (associative) monoidal $\infty$-category.
The underlying $A$-linear $\infty$-category $\LMod_B$
is equivalent to $\RMod_B=\LMod_{B^{op}}$ as an object of $\PR_A$.
Note that the underlying $\eone$-algebra is given by the restriction
along the map of $\infty$-operads $\eone^\otimes\to \etwo^\otimes$ determined by a linear injecitve map $\RR\to \RR^2$.
Any linear inclusion $\RR\to \RR^2$ is homotopy equivalent to one another via rotations.
Thus, we have a (non-canonical) equivalence $\LMod_B\simeq \RMod_B$ in $\PR_A$.

We regard $\LMod_B^\otimes$ as an object of $\Alg_1(\PR_A)$
and consider the opposite algebra $\LMod_B^{\otimes op}\in \Alg_1(\PR_A)$.
Namely, $\LMod_B^{\otimes op}$ has the ``opposite monoidal structure".
Let $\RMod_B^{\otimes}$ denote the monoidal $A$-linear presentable $\infty$-category
defined as $\Alg_1(\Theta_A)\circ\Alg_1(\iota_A)(B)\in \Alg_1(\PR_A)$.

\begin{Proposition}
\label{oppositemonoidal}
There exists an equivalence
\[
\LMod_B^{\otimes op}\simeq \RMod_B^\otimes
\]
in $\Alg_1(\PR_A)$.
\end{Proposition}

\Proof
We consider the diagram
\[
\xymatrix{
\Alg_1(\Alg_1(\Mod_A))  \ar[r]^(0.6){\Alg_1(\Theta_A)} \ar@<0.3ex>[d]_{\iota^{\textup{out}}}  \ar@<-0.3ex>[d]^{\iota^{\textup{in}}}  &  \Alg_1(\PR_A) \ar[d]^{\iota_{\PR_A}} \\
\Alg_1(\Alg_1(\Mod_A))   \ar[r]^(0.6){\Alg_1(\Theta_A)}   &  \Alg_1(\PR_A)
}
\]
where $\iota_{\PR_A}$ is the equvalence determined by the opposite algebras in $\Alg_1(\PR_A)$.
The arrow $\iota^{\textup{in}}$ is the equivalence given
by taking the opposite algebras with respect to the inner $\Alg_1$.
The arrow $\iota^{\textup{out}}$ is the equivalence given
by taking the opposite algebras with respect to the outer $\Alg_1$.
The equivalence $\iota^{\textup{in}}$ is induced by an automorphism $a:\etwo^\otimes\to \etwo^\otimes$
determined by reflecting across a line on the disc.
Also, $\iota^{\textup{out}}$ is induced by an automorphism $b:\etwo^\otimes\to \etwo^\otimes$
determined by reflecting across a line on the disc. The difference is the choice of a line on the disc.
Thus, such reflections are homotopy equivalent to each other. 
Thus,  $\iota^{\textup{in}}$ is equivalent to  $\iota^{\textup{out}}$.
Since $\Theta_A$ is a symmetric monoidal functor, there exists a canonical equivalence 
$\Alg_1(\Theta_A)\circ \iota^{\textup{out}}\simeq \iota_{\PR_A} \circ \Alg_1(\Theta_A)$.
By definition, $\RMod_B^\otimes=\Alg_1(\Theta_A)\circ\iota^{\textup{in}}(B)$
and $\LMod_B^{\otimes op}=\iota_{\PR_A}\circ \Alg_1(\Theta_A)(B)$.
It follows that $\RMod_B^\otimes\simeq \LMod_B^{\otimes op}$.
\QED

\begin{Construction}
\label{thetaTTT}
For a later use, we will construct
$\Theta_k^{+}:\Alg_{2}(\Mod_k)\to  \Alg_{1}(\Mod_{\Mod_k}(\PRTT))$.
which carries $B$ to $(\LMod_B^\otimes, \LMod_B^{\le0})$,
where $\LMod_B^{\le0}$ is the smallest full subcategory which contains $B$ and is closed under small colimits and extensions (see Remark~\ref{unifiedtstr}).
Recall that
$\Theta_k:\Alg_{2}(\Mod_k)\to  \Alg_{1}(\Mod_{\Mod_k}(\PR))$
is fully faithful (cf. Section~ref{MOA}). We denote the essential image by $T$. Let $S$ be the full subcategory of $\Alg_{1}(\Mod_{\Mod_k}(\PRTT))$
spanned by objects of the form $(\LMod_B^\otimes, \LMod_B^{\le0})$ with $B\in\Alg_{2}(\Mod_k)$.
The forgetful functor $\Alg_{1}(\Mod_{\Mod_k}(\PRTT))\to \Alg_1(\PR_k)$ induces $f:S\to T$.
We claim that $f$ is an equivalence. 
We define $\Theta_k^{+}$ to be the composition $g\circ \Theta_k|_{\Alg_{2}(\Mod_k^{\le0})}$
where $g$ is an inverse functor of $f$.
Clearly, $f$ is essentially surjective.
Note that every colimit-preserving monoidal functor $F:\LMod_A\to \LMod_B$
sends $A$ to $B$. Since $\LMod_A^{\le0}$ (resp. $\LMod_B^{\le0}$)
is the smallest full subcategory of $\LMod_A$ (resp. $\LMod_B$)
which contains $A$ (resp. $B$) and is closed under small colimits and extensions,
it follows that $F$ is right $t$-exact. Thus, $f$ is fully faithful.

We define
$\Theta_k^{\pm}:\Alg_{2}(\Mod_k^{\le0})\to  \Alg_{1}(\Mod_{\Mod_k}(\PRTTT))$
to be the restriction of $\Theta_k^+$.
If $B$ is connective, the $t$-strucutre on $\LMod_B$ is left complete
(Example~\ref{connectivetstr}) so that the target is $\Alg_{1}(\Mod_{\Mod_k}(\PRTTT))$.
\end{Construction}

{\it Endomorphism algebras.}
Let $\mathcal{M}$ be an $\infty$-category.
Let $\mathcal{A}^\otimes$ be a monoidal $\infty$-category.
Suppose that $\MMM$ is a left $\mathcal{A}^\otimes$-module, that is,
it has a left $\mathcal{A}^\otimes$-module action in the symmetric monoidal $\infty$-category $\wCat$ where the monoidal $\infty$-category $\mathcal{A}^\otimes$ is regarded as an associative algebra
object in $\wCat$.
Let $M$ be an object of $\MMM$.
The endomorphism algebra object of $M$ in $\mathcal{A}$ is defined as an associative
algebra $\End^{\mathcal{A}}(M)$ in $\mathcal{A}$ (i.e., an object of $\Alg_1(\mathcal{A})$)
together with a left module action of $\End^{\mathcal{A}}(M)$ on $M$, 
which has the universal property described informally as follows:
For any associative algebra $R$ in $\mathcal{A}$ and a left module action $R\otimes M\to M$
on $M$, there exists an essentially unique morphism $f:R\to \End^{\mathcal{A}}(M)$
such that the left module action of $R$ is induced by the restriction of the action of $\End^{\mathcal{A}}(M)$ along $f$ (this is a rough account, see \cite[4.7.1]{HA} for the precise formulation).

\subsection{Linear stable $\infty$-categories.}
\begin{Notation}
Write $(\PRTT)_k=\Mod_{\Mod_k}(\PRTT)$ and  $(\PRTTT)_k=\Mod_{\Mod_k}(\PRTTT)$.
\end{Notation}

Let $B\in \Alg_2(\Mod_k)$.
We define several versions of $B$-linear stable $\infty$-categories.

\begin{Definition}
\begin{enumerate}
\item
Let $\LMod^\otimes_B\in \Alg_1(\PR_k)$ be the associative monoidal $\infty$-category of left $B$-modules in $\Mod_k$.
Let $\RMod_{\LMod_B^\otimes}(\PR_k)$ be
the $\infty$-category of right $\LMod_B^\otimes$-modules in $\PR_k$.
We refer to an object of $\RMod_{\LMod_B^\otimes}(\PR_k)$
as a right $B$-linear stable $\infty$-category.

\item 
Suppose that $B$ is coconnective.
We consider the $t$-structure $(\LMod_B^{\le0}, \LMod_B^{\ge0})$ given in Proposition~\ref{tstrcoconnective}.
Note that $\LMod_B^{\le0}$ is closed under tensor products, and the unit object $B$ belongs to $\LMod_B^{\le0}$.
Consequently, $(\LMod_B,\LMod_B^{\le0})$ is naturally promoted to an object of $\Alg_1(\PRTT)$.
Let $\LMod_{\LMod_B^\otimes}((\PRTT)_k)$ denote
the $\infty$-category of right  $(\LMod_B,\LMod_B^{\le0})$-modules in $(\PRTT)_k$.
We refer to an object of $\LMod_{\LMod_B^\otimes}((\PRTT)_k)$
as a left $B$-linear stable $\infty$-category with $t$-stucture.

\item 
Suppose that $B$ is connective.
We consider the $t$-structure $(\LMod_B^{\le0}, \LMod_B^{\ge0})$ given in Example~\ref{connectivetstr}.
Note that $\LMod_B^{\le0}$ is closed under tensor products, and the unit object $B$ belongs to $\LMod_B^{\le0}$.
Consequently, $(\LMod_B,\LMod_B^{\le0})$ is naturally promoted to an object of $\Alg_1(\PRTTT)$.
Let $\RMod_{\LMod_B^\otimes}((\PRTTT)_k)$ denote
the $\infty$-category of right $(\LMod_B,\LMod_B^{\le0})$-modules in $(\PRTTT)_k$.
We refer to an object of $\RMod_{\LMod_B^\otimes}((\PRTTT)_k)$
as a complete right $B$-linear stable $\infty$-category.

\end{enumerate}

\end{Definition}

\begin{Remark} 
The left completion functor $\HL$ induces
\[
\HL_B:\RMod_{\LMod_B^\otimes}((\PRTT)_k) \to \RMod_{\HL(\LMod_B^\otimes)}((\PRTTT)_k)
\]
If $B$ is connective and $\LMod_B$ is
endowed with the left complete right complete $t$-structure in Example~\ref{connectivetstr},
then $\HL(\LMod_B^\otimes)\simeq \LMod_B^\otimes$.
\end{Remark}

\subsection{Koszul duals of algebras}
\label{koszulsec}
Let $\MMM^\otimes$ be a monoidal $\infty$-category. We write $\uni$ for a unit object.
We review the Koszul duality of algebras in $\MMM^\otimes$ (see \cite{HA}, \cite[X]{DAG} for more details).
Let $\Alg^+_{1}(\MMM)=\Alg_1(\MMM)_{/\uni}$ be the $\infty$-category of augmented $\eone$-algebra objects (associative algebra objects) in $\MMM^\otimes$.
For a pair $(\epsilon_A:A\to \uni,\epsilon_B:B \to \uni)\in \Alg^+_{1}(\MMM)\times \Alg^+_{1}(\MMM)$ we put
\[
\textup{Pair}(\epsilon_A,\epsilon_B):= \Map_{\Alg_{1}(\MMM)}(A\otimes B,\uni)\times_{\Map(A,\uni)\times\Map(B,\uni)}\{(\epsilon_A,\epsilon_B)\}.
\]
Fix $\epsilon_B:B \to \uni$ and consider the functor
$u_{\epsilon_B}:\Alg_1^+(\MMM)^{op}\to \wSSS$ given informally by $[\epsilon_A:A\to \uni]\mapsto \textup{Pair}(\epsilon_A,\epsilon_B)$. 
The Koszul dual $\DD_{\MMM^\otimes}(B)\in \Alg_1^+(\MMM)$ of $B$ (we abuse notation since $\DD_{\MMM^\otimes}(B)$ depends on the augmentation $\epsilon_B:B \to \uni$)
is defined by the universal property that $\DD_{\MMM^\otimes}(B)$ represents the functor $u_{\epsilon_B}$.
Under a mild condition on $\MMM^\otimes$, there exist $\DD_{\MMM^\otimes}(B)$ and functorial equivalences
\[
\Map_{\Alg_1^+(\MMM)}(A,\DD_{\MMM^\otimes}(B))\simeq \textup{Pair}(\epsilon_A,\epsilon_B)
\]
for $\epsilon_A:A\to \uni$.

Let $n\ge1$ be a natural number and  let $\Alg_{n}(\Mod_k)$ be the $\infty$-category of $\eenu$-algebras in $\Mod_k$.
We write $\Alg^+_{n}(\Mod_k)=\Alg_n(\Mod_k)_{/k}$ for the $\infty$-category
$\Alg_{n}(\Mod_k)_{/k}$ of augmented $\eenu$-algebras.
We consider the Koszul duals of augmented $\eenu$-algebras.
We take $\MMM^\otimes$ to be $\Alg_{n-1}^\otimes(\Mod_k)$. There exists the Koszul duality functor
\[
\DD_n:\Alg_{n}^+(\Mod_k)\simeq \Alg_1^+(\Alg_{n-1}(\Mod_k))\longrightarrow\Alg_1^+(\Alg_{n-1}(\Mod_k))\simeq  \Alg_{n}^+(\Mod_k)^{op}
\]
(see e.g. \cite[X, 4.4]{DAG}).
Let $\epsilon_{B}:B\to k$ be an augmented $\eenu$-algebra.
As mentioned above,
$\DD_n(B)$ is determined by the universal property.
For an augmented $\eenu$-algebra $\epsilon_C:C\to k$, we write
\[
\textup{Pair}(\epsilon_B,\epsilon_C)= \Map_{\Alg_{n}(\Mod_k)}(B\otimes_k C,k)\times_{\Map(B,k)\times\Map(C,k)}\{(\epsilon_B,\epsilon_C)\}.
\]
Then the functor $\Alg^+_{n}(\Mod_k)^{op}\to \wSSS$ given by $[\epsilon_C:C\to k]\mapsto \textup{Pair}(\epsilon_B,\epsilon_C)$
is representable by $\DD_{n}(B)$.
There is a universal ``Koszul dual pairing'' $B\otimes_k\DD_{n}(B)\to k$ corresponding to
the identity map $\textup{id}\in \Map_{\Alg^+_{n}(\Mod_k)}(\DD_{n}(B),\DD_{n}(B))$.
The Koszul dual $\DD_{n}(B)$ can also be interpreted as a centralizer of $B\to k$, see \cite[5.3.1]{HA}.
We refer to a comprehensive account \cite{Gin} for centralizers and  Koszul duals.
We refer to $\DD_{n}(B)$ as the $\eenu$-Koszul dual of $B$.

The Koszul duals can be described in terms of bar constructions.
Let $\textup{Bar}:\Alg_{1}^+(\Mod_k)\to \Mod_{k}$ be the functor
given by $[B\to k] \mapsto k\otimes_Bk$, which we refer to as the bar construction of augmented algebras.
Then $k\otimes_Bk$ admits a structure of a coaugmented coalgebra in a suitable way, and $\textup{Bar}$
is promoted to $\textup{Bar}:\Alg_{1}^+(\Mod_k)\to \Alg^+_{1}((\Mod_k)^{op})^{op}$
(notice that the target is the $\infty$-category of coaugmented coalgebras).
Interating bar construction we have
induces $\textup{Bar}^n:\Alg_{n}^+(\Mod_k)\to \Alg_{n}^+((\Mod_k)^{op})^{op}$.
According to \cite[X, 4.4.20]{DAG}, for $B\in \Alg^+_{n}(\Mod_k)$, the $\eenu$-Koszul dual $\DD_{n}(B)$
is equivalent to
the $k$-linear dual $\mathcal{H}om_{k}(\textup{Bar}^n(B),k)$.
In particular, $\eone$-Koszul dual of $B\to k$ is given by the endomorphism algebra $\End_B(k)$ with the canonical
augmentation $\End_B(k)\to \End_k(k)\simeq k$ (here $k$ is regarded as a left $B$-module).

\section{Koszul duality and left completion}

\label{adjointsec}

\subsection{}

Let
$A$ be an augmented $\etwo$-algebra in $\Mod_k$. Suppose that $A$ is connective.
Note that $\LMod_A^{\le0}$ is the smallest full subcategory which contains $A$ and is closed under small colimits (see Exmaple~\ref{connectivetstr}, \cite[7.1.1.13]{HA}).
Thus, $\LMod_A^{\le0}$ contains an unit object $R$ and is closed under tensor products
so that the monoidal $\infty$-category $\LMod_A^\otimes$
can be regarded as an algebra object in $\PRTTT$.
Let $\DDD$ be a complete $A$-linear stable $\infty$-category endowed with $(\DDD^{\le0},\DDD^{\ge0})$,
that is, $\DDD\in \RMod_{\LMod_A}(\PRTTT)$.
In particular, the $t$-structure is accessible, left complete, and right complete.
The tensor product of $\PRTT$ induces the right complete accessible $t$-structure on $\DDD\otimes_{\LMod_{A}}\Mod_k$.
This $t$-structure is determined by
$(\DDD\otimes_{\LMod_A}\Mod_k)^{\le0}$ defined as
the smallest full subcategory which contains a collection of objects $\{D\otimes_Ak\}_{D\in \DDD^{\le0}}$
and is closed under small colimits and extensions.
Here $D\otimes_Ak$ is the image of $(D,k)\in \DDD\times \Mod_k$ under the canonical functor $\DDD\times \Mod_k\to \DDD\otimes_{\LMod_A}\Mod_k$.

From \cite[4.8.4.6]{HA} there exists a canonical equivalence $\DDD\otimes_{\LMod_A}\Mod_k\simeq\RMod_k(\DDD)$. 
We consider the forgetful functor 
\[
u:\DDD\otimes_{\LMod_A}\Mod_k\simeq\RMod_k(\DDD) \to \DDD.
\]

\begin{Lemma}
\label{texactforget}
The functor $u:\DDD\otimes_{\LMod_A}\Mod_k\simeq\RMod_k(\DDD) \to \DDD$
is $t$-exact.
Moreover, $M$ lies in $(\DDD\otimes_{\LMod_A}\Mod_k)^{\le0}$
if and only if its image in $\DDD$ lies in $\DDD^{\le0}$.
\end{Lemma}

\Proof
Note that the forgetful functor $u$ preserves small colimits.
Note also that $D\otimes_Ak\in \DDD^{\le0}$ for any $D\in \DDD^{\le0}$
since $k$ lies in $\LMod_A^{\le0}$ and by definition the tensor product induces $\DDD^{\le 0}\times \LMod_A^{\le0} \to \DDD^{\le0}$.
Taking account of the definition of $(\DDD\otimes_{\LMod_A}\Mod_k)^{\le0}$
we conclude that $u$ is right $t$-exact.

Next, we will prove that $u$ is left $t$-exact.
Let $v:\DDD\to \DDD\otimes_{\LMod_A}\Mod_k$
be the left adjoint induced by the base change functor $\LMod_A\to \Mod_k$.
This left adjoint carries $D$ to $D\otimes_Rk$.
In particular, $\DDD^{\le0}$ maps to $(\DDD\otimes_{\LMod_A}\Mod_k)^{\le0}$.
It follows that $v$ is right $t$-exact and $u$ is left $t$-exact.
Moreover, $E\in\DDD\otimes_{\LMod_A}\Mod_k $
 lies in $(\DDD\otimes_{\LMod_A}\Mod_k)^{\ge0}$
if and only if $\Hom_{\hhh(\DDD\otimes_{\LMod_A}\Mod_k)}(D\otimes_Ak[1],E)$
for any $D\in \DDD^{\le0}$.
Thus, $(\DDD\otimes_{\LMod_A}\Mod_k)^{\ge0}=u^{-1}(\DDD^{\ge0})$.
\QED

Strictly speaking, the following is not necessary for Theorem~\ref{intromain1}, but useful.

\begin{Proposition}
\label{already}
Let $\DDD$ be an object of $\RMod_{\LMod_A}(\PRTTT)$ as above.
Let $\DDD\otimes_{\LMod_A}\Mod_k\in (\PRTT)_k$
denote the relative tensor product in $\PRTT$.
Equivalently, $\DDD\otimes_{\LMod_A}\Mod_k$
is the relative tensor product in $\PR_{\mathbb{S}}$
endowed with the $t$-structure determined by $(\DDD\otimes_{\LMod_A}\Mod_k)^{\le0}$
defined above.
Then the $t$-structure $(\DDD\otimes_{\LMod_A}\Mod_k, (\DDD\otimes_{\LMod_A}\Mod_k)^{\le0})$
is left complete.
\end{Proposition}

\Proof
For $M\in \DDD\otimes_{\LMod_A}\Mod_k$ we consider the natural morphism
$\theta:M\to \lim_{n\to -\infty}(\tau^{\ge n}M)$.
It suffices to prove that $\theta$ is an equivalence.
Note that $u:\DDD\otimes_{\LMod_R}\Mod_k\to \DDD$ is $t$-exact and preserves small colimits (see Lemma~\ref{texactforget}),
and $\DDD$ and $\DDD\otimes_{\LMod_A}\Mod_k$ are right complete.
Thus we have
\begin{eqnarray*}
u(M) &\to& u(\lim_{n\to -\infty}(\tau^{\ge n}M)) \\
&\simeq& \lim_{n\to -\infty}u((\tau^{\ge n}M)) \\
&\simeq& \lim_{n\to -\infty}\tau^{\ge n}u(M).
\end{eqnarray*}
By the left completeness of $(\DDD,\DDD^{\le0})$,
the composite is an equivalence.
Thus, $u(M) \to u(\lim_{n\to -\infty}(\tau^{\ge n}M))$
is an equivalence.
The forgetful functor $u$ is conservative so that $\theta$
is an equivalence.
\QED

\begin{Remark}
\label{alreadyrem}
In Lemma~\ref{texactforget}, $M$ lies in $(\DDD\otimes_{\LMod_A}\Mod_k)^{\ge0}$
if and only if its image in $\DDD^{\ge0}$.
To see this,
we write the canonical morphism $p:\tau^{\le0}M\to M$
from the truncation (the connective cover).
The morphism $p$ is an equivalence if and only if $M\in (\DDD\otimes_{\LMod_A}\Mod_k)^{\le0}$.
By the $t$-exactness of $u$ and the left completeness of the $t$-structures on $\DDD$ and $\DDD\otimes_{\LMod_A}\Mod_k$ (Proposition~\ref{already}),
we have $u(p):u(\tau^{\le0}M)\simeq\tau^{\le0}u(M)\to u(M)$.
The morphism $u(p)$ is an equivalence if and only if $u(M)\in \DDD^{\le0}$.
Since $u$ is conservative, $M$ lies in $(\DDD\otimes_{\LMod_R}\Mod_k)^{\le0}$ exactly when $u(M)\in \DDD^{\le0}$.
Consequently, $M$ lies in $(\DDD\otimes_{\LMod_R}\Mod_k)^{\le0}$ (resp. $(\DDD\otimes_{\LMod_R}\Mod_k)^{\ge0}$)
if and only if its image in $\DDD$ lies in $\DDD^{\le0}$ (resp. $\DDD^{\ge0}$).
\end{Remark}

\begin{Lemma}
\label{tensorvshom}
Let $\DDD$ be a right $\LMod_A$-module in $\PR_{\mathbb{S}}$, that is an object of $\RMod_{\LMod_A}(\PR)$.
There exists a canonical equivalence
\[
\DDD\otimes_{\LMod_{A}}\Mod_k\simeq \Fun_{\LMod_A}^r(\Mod_k, \DDD)
\]
in $\PR_{\mathbb{S}}$. Here the right-hand side is the functor category from $\Mod_k$ to $\DDD$ in  $\RMod_{\LMod_A}(\PR)$.
\end{Lemma}

\Proof
According to \cite[4.8.4.1, 4.8.4.6]{HA}, we have equivalences $\DDD\otimes_{\LMod_{A}}\Mod_k\simeq \RMod_k(\DDD)\simeq \Fun_{\LMod_A}^r(\Mod_k, \DDD)$.
\QED

\begin{Remark}
Let $\DDD \in \RMod_{\LMod_A}(\PRTTT)$. 
It is natural to consider a relationship between the hearts $\DDD^{\heartsuit}$ and $(\DDD\otimes_{\LMod_A}\Mod_k)^\heartsuit$.
Observe that $u:\DDD\otimes_{\LMod_A}\Mod_k\to \DDD$ is equivalent to $\DDD\otimes_{\LMod_A}\Mod_k\simeq \Fun_{\LMod_A}^r(\Mod_k, \DDD)\to \DDD$
given by the evaluation at $k\in \Mod_k$.
According to Lemma~\ref{texactforget} and Remark~\ref{alreadyrem}, the heart
$\Fun_{\LMod_A}^r(\Mod_k, \DDD)^{\heartsuit}$ (induced by the $t$-structure on $\DDD\otimes_{\LMod_{A}}\Mod_k$)
consists of those functors $F$ such that $F(k)\in \DDD^{\heartsuit}$.
Put another way, $(\DDD\otimes_{\LMod_a}\Mod_k)^\heartsuit\simeq \RMod_k(\DDD)\times_{\DDD}\DDD^\heartsuit$
where $\RMod_k(\DDD)\to\DDD$
is the forgetful functor.
The restriction to the heart $\Mod_k^{\heartsuit}$ (that is, the abelian category of $k$-vector spaces)
determines the exact functor of abelian categories
\[
\Gamma:\Fun_{\LMod_A}^r(\Mod_k, \DDD)^{\heartsuit}\to \textup{Lin}_{H^0(A)}(\Mod_k^\heartsuit, \DDD^{\heartsuit})
\]
where the right-hand side is the abelian category of 
$H^0(A)$-linear coproduct-preserving exact functors.
One can prove that 
the functor $\Gamma$ is an equivalence.
In addition, one easily sees that they are canonically equivalent to the full subcategory of $\DDD^\heartsuit$ spanned by objects $D$ such that the $H^0(A)$-linear structure morphism $H^0(A)\to \End_{\DDD^\heartsuit}(D)$
factors through $H^0(A)\to k$. Here $\End_{\DDD^\heartsuit}(D)$ is the ordinary endomorphism algebra.
As a consequence one may identify the heart $(\DDD\otimes_{\LMod_A}\Mod_k)^\heartsuit$ with 
$\textup{Lin}_{H^0(A)}(\Mod_k^\heartsuit, \DDD^{\heartsuit})$.
The latter functor category can be interpreted as the $!$-pullback
of $\DDD^{\heartsuit}$ along $\Spec k\to \Spec H^{0}(A)$.
We 
regard $\DDD$
as a deformation of $\EE\simeq \DDD\otimes_{\LMod_A}\Mod_k\in (\PRTTT)_k$ to $R$.
The heart $\DDD^{\heartsuit}$ can be regarded as a deformation of the $k$-linear abelian 
category $\EE^\heartsuit$ to a $H^0(A)$-linear abelian category.
We will study the deformations of the heart induced by deformations of $\DDD$ elsewhere.

\end{Remark}

\subsection{}
Let $A$ be an augmented $\etwo$-algebra in $\Mod_k$. 
We write $\epsilon_A:A\to k$ for the augmentation.
We think of $\LMod_A$ as a monoidal $\infty$-category (cf. Section~\ref{MOA}).
Consider $\Mod_k$ to be the left $\LMod_A$-module object via the monoidal functor $\LMod_A^\otimes\to \LMod_k^\otimes=\Mod_k^\otimes$.
Let $\End^l_{A}(\Mod_k)$ be the endomorphism algebra object of $\Mod_k$ in $\LMod_{\LMod_A}(\PR)$.
This is defined as an object of $\Alg_1(\PR_k)$, which has a universal left module action on $\Mod_k$.
When we emphasize the monoidal structure, we write $\End^l_{A}(\Mod_k)^\otimes$ for it.
There exists a $\eone$-Koszul dual pairing
\[
\kappa:\LMod_A^\otimes\otimes_k\End^l_{A}(\Mod_k)^\otimes\longrightarrow \Mod_k^\otimes
\]
which is a morphism of $\Alg_1(\PR_k)$.

\begin{Lemma}
\label{dualcatalg}
The $\etwo$-Koszul dual $\DD_2(A)$ can be identified with $E_k(\End^l_{A}(\Mod_k))\in \Alg_2(\Mod_k)$ (see Section~\ref{MOA} for $E_k$).
The universal $\etwo$-Koszul pairing $A\otimes_k\DD_2(A) \to k$
is equivalent to the composition of $E_k(\kappa)$ and the canonical 
morphism $E_k(\LMod_A)\otimes_k E_k(\End^l_{A}(\Mod_k))\to E_k(\LMod_A\otimes_k\End^l_{A}(\Mod_k))$.
\end{Lemma}

\Proof
We first recall the universal property of $\DD_2(A)$. 
For $\epsilon_A:A\to k$ and $\epsilon_B:B\to k$ we consider the space
\[
\textup{Pair}(\epsilon_A,\epsilon_B):= \Map(A\otimes_k B,k)\times_{\Map(A,k)\times\Map(B,k)}\{(\epsilon_A,\epsilon_B)\}.
\]
Here $\Map$ means the mapping space in $\Alg_{2}(\Mod_k)$.
Consider the functor
$u_{\epsilon_A}:\Alg_2^+(\Mod_k)^{op}\to \wSSS$ given informally by $[\epsilon_B:B\to k]\mapsto \textup{Pair}(\epsilon_A,\epsilon_B)$. 
The Koszul dual $\DD_{2}(A)\in \Alg_2^+(\Mod_k)$ of $A$
is defined by the universal property that $\DD_{2}(A)$ represents the functor $u_{\epsilon_A}$.
Consider $\LMod(\epsilon_A):\LMod_A\to \Mod_k$ 
and $\LMod(\epsilon_B):\LMod_B\to \Mod_k$.
Let $\epsilon_{\MMM^\otimes}:\MMM^\otimes\to \Mod_k^\otimes$ be an object of $\Alg_1^+(\PR_k)=\Alg_1(\PR_k)_{/\Mod^\otimes_k}$.
We define $\textup{Pair}(\LMod(\epsilon_A),\epsilon_{\MMM^\otimes})$ to be
\[
 \Map(\LMod_A^\otimes\otimes \MMM^\otimes,\Mod_k^\otimes)\times_{\Map(\LMod_A^\otimes,\Mod_k^\otimes)\times\Map(\MMM^\otimes,\Mod_k^\otimes)}\{(\LMod(\epsilon_A),\epsilon_{\MMM^\otimes})\}.
\]
Here $\Map$ means the mapping space in $\Alg_{1}(\PR_k)$.
Let 
$u_{\LMod(\epsilon_A)}:\Alg_1^+(\PR_k)^{op}\to \wSSS$
be the functor given informally by $[\epsilon_{\MMM^\otimes}:\MMM^\otimes\to \Mod_k^\otimes]\mapsto \textup{Pair}(\LMod(\epsilon_A),\epsilon_{\MMM^\otimes})$.
This functor
is represented by the $\eone$-Koszul dual
$\End^l_{A}(\Mod_k)$ with the canonical augmentation $\End^l_{A}(\Mod_k)\to\Mod_k$.

Since $\Theta_k:\Alg_2(\Mod_k)\to \Alg_1(\PR_k)$ is the fully faithful symmetric monoidal functor,
we have a canonical equivalence
$\textup{Pair}(\epsilon_A,\epsilon_B)\simeq \textup{Pair}(\LMod(\epsilon_A),\LMod(\epsilon_B))$.
Namely, $u_{\epsilon_A}$ is naturally equivalent to
the composite $u_{\LMod(\epsilon_A)}\circ \Theta_k^{op}:\Alg_2^+(\Mod_k)^{op}\to \Alg_1^+(\PR_k)^{op}\to \wSSS$.
Take $D=E_k(\End^l_{A}(\Mod_k))$ and the composite
\[
P:E_k(\LMod_A)\otimes_k E_k(\End^l_{A}(\Mod_k))\to E_k(\LMod_A\otimes_k\End^l_{A}(\Mod_k))\to E_k(\Mod_k)\simeq \Mod_k
\]
where the left arrow is induced by the lax symmetric monoidal functor $E_k$,
and the right arrow is induced by the universal pairing $\kappa$.
Note also that $E_k(\LMod_A)\simeq \LMod_A$.
Let $p:A\otimes_kD\to k$ be the morphism in $\Alg_2(\Mod_k)$
that corresponds to $P$ through $\Theta_k$.
Then taking account of the adjunction $\Theta_k:\Alg_2(\Mod_k)\rightleftarrows \Alg_1(\PR_k):E_k$, we see that $D$ together with $p$ defines a $\etwo$-Koszul dual of $A\to k$.
\QED

Using the counit of the adjunction we have:

\begin{Corollary}
\label{monoidal}
There exists a morphism $\xi:\LMod_{\DD_2(A)}^\otimes\to \End_{A}^{l}(\Mod_k)^\otimes$ in $\Alg_1(\PR_k)$.
\end{Corollary}

\begin{Lemma}
\label{moduleproduct}
Let $A$ be an augmented $\etwo$-algebra in $\Mod_k$.
There exists an equivalence
\[
\Mod_k\otimes_{\RMod_A}\Mod_k \simeq \RMod_{k\otimes_Ak}
\]
in $\PR_k$.
\end{Lemma}

\Proof
The functor $\Theta_k:\Alg_1(\Mod_k)\to (\PR_k)_{\Mod_k/}$
is a (fully faithful) symmetric monoidal left adjoint functor, and 
the forgetful functor $(\PR_k)_{\Mod_k/}\to \PR_k$ preserves geometric realizations.
Recall that the relative tensor product is
defined as the colimit of the simplicial diagram of the bar construction (cf. \cite{HA}).
Thus, our assertion follows.
\QED

\begin{Remark}
\label{moduleproduct2}
There also exists an equivalence 
\[
\Mod_k\otimes_{\RMod_A}\Mod_k\simeq \Mod_k\otimes_{\LMod_A}\Mod_k\simeq \Fun_{\LMod_A}^l(\Mod_k, \Mod_k)
\]
(see Proposition~\ref{tensorvshom} and Proposition~\ref{oppositemonoidal}).
\end{Remark}

\subsection{}
Suppose that
$A$ is connective and equipped with an augmentation $A\to k$.
Let $\DD_2(A)\in \Alg_2^+(\Mod_k)$ be the $\etwo$-Koszul dual of $A\to k$ (cf. Section~\ref{koszulsec}).
Since $A$ is connective, $\DD_2(A)$
is coconnective. Thus, 
$\LMod_{\DD_2(A)}$ has the $t$-structure $(\LMod_{\DD_2(A)}^{\le0}, \LMod_{\DD_2(A)}^{\ge0})$
defined in Proposition~\ref{tstrcoconnective}.
We regard $\LMod_{\DD_2(A)}$ as an associative algebra object in $\PRTT$.

Note that $k\otimes_Ak$ is also connective and $H^0(k\otimes_Ak)=k$.
Then
$\RMod_{k\otimes_Ak}$ has the standard $t$-structure, given in Example~\ref{connectivetstr},
such that the heart $\RMod_{k\otimes_Ak}^{\heartsuit}$
is the abelian category of $k$-vector spaces.
Let $\textup{RCoh}(k\otimes_Ak)$ be the smallest stable subcategory
which contains $k\in \RMod_{k\otimes_Ak}$ and is closed under retracts.
An object of $\textup{RCoh}(k\otimes_Ak)$
is characterized as an object $M$ such that 
it has bounded amplitude (i.e., $H^n(M)=0$ for $|n|>>0$)
and each $H^n(M)$ is finite dimensional (i.e. $\dim H^n(M)<<\infty$ for $n\in \ZZ$).

\begin{Lemma}
\label{leftcomplemmma}
Suppose that $A$ is connective.
The functor 
$\xi:\LMod_{\DD_2(A)}\to  \End_{A}^{l}(\Mod_k)$ in Corollary~\ref{monoidal} is a left completion.
\end{Lemma}

\Proof
By Lemma~\ref{moduleproduct} and Remark~\ref{moduleproduct2},
there exists an equivalence
$\End_{A}^{l}(\Mod_k)\simeq \RMod_{k\otimes_Ak}$ in $\PR_k$.
Note that $k\in \textup{RCoh}(k\otimes_Ak)$ is a single compact generator of the 
Ind-category $\Ind(\textup{RCoh}(k\otimes_Ak))$.
Consequently, if we denote by $\DD_1(k\otimes_Ak)$
the $\eone$-Koszul dual, 
there exists a canonical Morita-type equivalence $\textup{RCoh}(k\otimes_Ak)\simeq \textup{LPerf}_{\DD_1(k\otimes_Ak)}$
induced by the integral kernel $(k\otimes_Ak)^{op}$-$\DD_1(k\otimes_Ak)$-bimodule $k$,
where $\textup{LPerf}_{\DD_1(k\otimes_Ak)}$ is the smallest stable subcategory of $\textup{LMod}_{\DD_1(k\otimes_Ak)}$ which contains $\DD_1(k\otimes_Ak)$
and is closed under retracts.
Explicitly, the equivalence $\textup{LPerf}_{\DD_1(k\otimes_Ak)}\to \textup{RCoh}(k\otimes_Ak)$
is given informally by $N\mapsto k\otimes_{\DD_1(k\otimes_Ak)}N$.
We then have equivalences
\[
\Ind(\textup{RCoh}(k\otimes_Ak)) \simeq \Ind(\textup{LPerf}_{\DD_1(k\otimes_Ak)})\simeq \LMod_{\DD_1(k\otimes_Ak)} 
\]
in $\PR_k$, which carries $k$ to $\DD_1(k\otimes_Ak)$.
Moreover, $\DD_2(A)\simeq \DD_1(\textup{Bar}(A))\simeq \DD_1(k\otimes_Ak)=\End_{k\otimes_Ak}(k)$.
Thus, $\Ind(\textup{RCoh}(k\otimes_Ak))\simeq  \LMod_{\DD_2(A)}$ (they are also equivalent to $\RMod_{\DD_2(A)}$ in $\PR_k$ since $\DD_2(A)$ is an $\etwo$-algebra).
The functor $\xi$
is the colimit preserving functor $F:\Ind(\textup{RCoh}(k\otimes_Ak))\to \RMod_{k\otimes_{A}k}$
which extends $\textup{RCoh}(k\otimes_Ak)\hookrightarrow \RMod_{k\otimes_{A}k}$.
Thus, it is enough to show that $F$
 is a left completion functor.
Note that $\Ind(\textup{RCoh}(k\otimes_Ak))^{\le0}$ induces by $\LMod_{\DD_2(R)}^{\le0}$
is the smallest full subcategory which contains $k$ and is closed under small colimits and extensions.
If $(\textup{RCoh}(k\otimes_Ak)^{\le0}, \textup{RCoh}(k\otimes_Ak)^{\ge0})$ is the $t$-structure
induced by $\RMod_{k\otimes_Ak}$, $(\Ind(\textup{RCoh}(k\otimes_Ak)^{\le0}), \Ind(\textup{RCoh}(k\otimes_Ak)^{\ge0}))$
defines a $t$-structure on $\Ind(\textup{RCoh}(k\otimes_Ak))$ such that $\Ind(\textup{RCoh}(k\otimes_Ak)^{\le0})$
is the full subcategory which
contains $k$ and is closed under small colimits and extensions.
Thus, $\Ind(\textup{RCoh}(k\otimes_Ak))^{\le0}=\Ind(\textup{RCoh}(k\otimes_Ak)^{\le0})$.
Since $\RMod_{k\otimes_Ak}^{\ge0}$ is closed under filtered colimits,
it follows from the definitions of $t$-structures that $F$ is $t$-exact.
Since the $t$-structure on $\RMod_{k\otimes_Ak}$ is left complete,
it is enough to  prove that 
$F'=F|_{\Ind(\textup{RCoh}(k\otimes_Ak))^{\ge0}}:\Ind(\textup{RCoh}(k\otimes_Ak))^{\ge0}\to \RMod_{k\otimes_Ak}^{\ge0}$
is an equivalence.
To this end, we follow the steps in \cite[Proposition 1.2.4]{Ind}.
We first prove that $F'$ is fully faithful. We will prove that
\[
\Map_{\Ind(\textup{RCoh}(k\otimes_Ak))}(M,N)\to \Map_{\RMod_{k\otimes_Ak}}(F(M),F(N))
\]
is an equivalence
for any $N\in \Ind(\textup{RCoh}(k\otimes_Ak))^{\ge0}$ and any $M\in \Ind(\textup{RCoh}(k\otimes_Ak))$.
We may and will assume that $M$ lies in $\textup{RCoh}(k\otimes_Ak)$.
Also, we may take $N$ to be a filtered colimit $\colim_{i\in I}N_i$ such that $N_i\in \textup{RCoh}(k\otimes_Ak)^{\ge0}$.
We claim that $\Map_{\RMod_{k\otimes_Ak}}(M,-)$ commutes with filtered colimits in $\RMod_{k\otimes_Ak}^{\ge0}$.
Our claim shows the desired equivalence
\begin{eqnarray*}
\Map_{\Ind(\textup{RCoh}(k\otimes_Ak))}(M,\colim_{i\in I}N_i) &\simeq& \colim_{i\in I}\Map_{\Ind(\textup{RCoh}(k\otimes_Ak))}(M,N_i) \\
&\simeq& \colim_{i\in I}\Map_{\RMod_{k\otimes_Ak}}(M,N_i)  \\
&\simeq& \Map_{\RMod_{k\otimes_Ak}}(M,\colim_{i\in I}N_i).
\end{eqnarray*}
To prove the claim, we suppose that the amplitude of $M$ is contained in $[-\infty,n]$.
When $n=-1$, the claim is clear because the mapping spaces are contractible.
We will assume that the claim holds when the amplitude of $M$ is contained in $[-\infty,n-1]$.
We will prove the case when the amplitude of $M$ is contained in $[m,n]$.
Let $r$ be the dimension of $H^n(M)$.
Take $\alpha:(k\otimes_Ak)^{\oplus r}[-n]\to M$
such that it induces an isomorphism $k^{\oplus r}=H^n((k\otimes_Ak)^{\oplus r}[-n])\to H^n(M)$.
It gives rise to the cofiber sequence
\[
(k\otimes_Ak)^{\oplus r}[-n]\to M  \to C
\]
such that 
the amplitude of the cofiber $C$ is contained in $[-\infty,n-1]$.
From the assumption and the compactness of $(k\otimes_Ak)^{\oplus r}[-n]$ in $\RMod_{k\otimes_Ak}$
the claim holds for the case when the amplitude of $M$ is contained in $[-\infty,n]$.
Consequently, we conclude that $F'$ is fully faithful.
Next, we will show that $F'$ is essentially surjective.
This immediately follows from the fact that every object of $\RMod_{k\otimes_Ak}^{\ge0}$
is a filtered colimits of $\colim_{i\in I}N_i$ with $N_i\in \textup{RCoh}(k\otimes_Ak)^{\ge0}$
(use the right completeness of the $t$-structure on $\RMod_{k\otimes_Ak}$ and the characterization of objects of $\textup{RCoh}(k\otimes_Ak)$).
\QED

\begin{Definition}
Let $R$ be an augmented $\eenu$-algebra in $\Mod_k$. Let $n$ be a natural number or $\infty$.
The augmented $\eenu$-algebra $R$ is Artin if the following properties hold:
\begin{itemize}
\item $R$ is connective,

\item $\oplus_{n\in \ZZ}H^{n}(R)$ is finite dimensional over $k$,

\item if $I$ is the radical of $H^0(R)$, the canonical map $k\to H^0(R)/I$ is an isomorphism.

\end{itemize}
By convention, an $\mathbf{E}_{\infty}$-algebra is a commutative algebra object in $\Mod_k$, and $\Alg_\infty(\Mod_k)=\CAlg(\Mod_k)$.
We let $\textup{Art}_n$ denote the full subcategory on $\Alg_n^+(\Mod_k)$ spanned by Artin $\eenu$-algebras.
\end{Definition}

In \cite{DAG}, the terminology small is adopted instead of Artin.

\begin{Lemma}
\label{dbar}
If $R$ is an Artin $\etwo$-algebra, there
exists an equivalence $k\otimes_{\DD_2(R)}k\simeq \DD_1(R)$ in $\Alg_1(\Mod_k)$.
\end{Lemma}

\Proof
Since $R$ is Artin, the biduality morphisms $R\to \DD_1\DD_1(R)$ and $R\to \DD_2\DD_2(R)$ are equivalences:
This is an immediate consequence of \cite[X, Theorem 4.4.5, Proposition 4.5.1,  Proposition 4.5.6, Propositions 4.1.14 and 4.1.13]{DAG}.
It also follows from these results that $\DD_2(R)$ is locally finite.

Note that $\DD_1(k\otimes_{\DD_2(R)}k)=\DD_1(\textup{Bar}(\DD_2(R))\simeq \DD_2\DD_2(R)\simeq R$.
Thus, $\DD_1\DD_1(R)\simeq R\simeq \DD_1(k\otimes_{\DD_2(R)}k)$.
Let $\Alg_1^+(\Mod_k)^{\textup{cclf}}$ be the full subcategory
of $\Alg_1^+(\Mod_k)$ spanned by coconnective and locally finite algebras.
According to \cite[X, 3.1.15]{DAG}, 
(the restriction of) the Koszul duality functor $\DD_1:\Alg_1^+(\Mod_k)^{\textup{cclf}}\to \Alg_1^+(\Mod_k)^{op}$
is fully faithful. 
For now, we suppose that $k\otimes_{\DD_2(R)}k$ is coconnective and locally finite.
It follows that $\DD_1(R)\simeq k\otimes_{\DD_2(R)}k$.

It remains to verify that $k\otimes_{\DD_2(R)}k$ is coconnective and locally finite.
By \cite[X, 4.5.9]{DAG}, $\DD_2(R)$ is $2$-coconnective.
In particular, $H^1(\DD_2(R))=0$.
We write $B=\DD_2(R)$.
We will construct the sequence
\[
B(0)\to B(1)\to \cdots \to B(m)\to \cdots
\]
in $\Alg_k(\Mod_k)_{/k}$ such that $B(0)=B$,
and $B(m)$ is a locally finite $(m+2)$-coconnective $B$-module which fits into the pushout square
\[
\xymatrix{
B^{\oplus d_{m}}[-m-2] \ar[r] \ar[d] & 0 \ar[d] \\
 B(m) \ar[r] &  B(m+1)
}
\]
in $\LMod_B$ with $d_{m}\ge0$.
Take $B(0)\to k$ to be the augmentation $B\to k$.
We suppose that a locally finite $(m+2)$-coconnective $B$-mdoule $B(m)$ with $B(m)\to k$ has been constructed.
Let $d_m$ be the dimension of $H^{m+2}(B(m))$.
Choose $f_m:B^{\oplus d_m}[-m-2]\to B(m)$
which induces an isomorphism $k^{\oplus d_m}=H^{m+2}(B^{\oplus d_m}[-m-2])\to H^{m+2}(B(m))$.
Let $B(m+1)$ be the cofiber of $f_m$ in $\LMod_B$, and let $B(m)\to B(m+1)$ be the canonical morphism.
Using the long exact sequence of cohomology groups,
we see that $B(m+1)$ is $(m+3)$-coconnective and locally finite.
Note that the filtered colimit $\colim_{m\to \infty}B(m)$ is canonically equivalent to $k$.
Thus, $k\otimes_Bk\simeq \colim_{m\to \infty} k\otimes_BB(m)$.
We will show that if $k\otimes_BB(m)$ is coconnective and locally finite, then $k\otimes_BB(m+1)$ is also coconnective and locally finite.
Note that $k\otimes_BB(m+1)$ is a cofiber of $k^{\oplus d_m}[-m-2]=k\otimes_BB^{\oplus d_{m}}[-m-2]\to k\otimes_BB(m)$.
The long exact sequence of cohomology groups shows that $k\otimes_BB(m+1)$ is also coconnective and locally finite.
Thus, $k\otimes_Bk$ is coconnective.
In addition, the argument shows that for any $i\in \ZZ$ the sequence 
\[
H^i(k\otimes_BB(0))\to \cdots \to H^i(k\otimes_BB(m))\to H^i(k\otimes_BB(m+1))\to \cdots
\]
is eventually constant: There exists a natural number $N$ such that $H^i(k\otimes_BB(m))\to H^i(k\otimes_BB(m+1))$
is an isomorphism for $m>N$. It follows that $k\otimes_Bk$ is locally finite.
\QED

\begin{Remark}
Since $R$ is an $\etwo$-algebra, $R\simeq R^{op}$ as $\eone$-algebras. In particular,
$\DD_1(R)\simeq \DD_1(R^{op})\simeq \DD_1(R)^{op}$.
\end{Remark}

\begin{Lemma}
\label{leftcomp2}
There exists a left completion functor $\eta:\RMod_{k\otimes_{\DD_2(R)}k}\to \LMod_R$.
Moreover, the target $\LMod_R$ is a free $\LMod_R$-module of rank one.
\end{Lemma}

\Proof
Let $\textup{LCoh}(R)\subset \LMod_R$ be the smallest stable subcategory which contains $k$
and is closed under retracts.
Since the heart $\LMod_R^\heartsuit$ is equivalent to the abelian category of $H^0(R)$-modules 
and $H^0(R)$ is an ordinary Artin local $k$-algebra with residue field $k$, $\textup{LCoh}(R)$
is the full subcategory spanned by
those objects $M$ such that 
$M$ has bounded amplitude
and each $H^n(M)$ is finite dimensional.
The presentbale $\infty$-category $\Ind(\textup{LCoh}(R))$
has a single compact generator $k$.
If we write $\End_R(k)$ for the endomorphism algebra of $k\in \Ind(\textup{LCoh}(R))$, then
$\End_R(k)\simeq\DD_1(R)$.
Thus, 
there exists an equivalence $\Perf_{\DD_1(R)}\simeq \textup{LCoh}(R)$
which carries $N$ to $N\otimes_{\DD_1(R)}k$ ($\Perf_{\DD_1(R)}$ denotes the smallest stable subcategory of $\RMod_{\DD_1(R)}$ which contains $\DD_1(R)$
and is closed under retracts).
Taking Ind-categories it gives rise to $\alpha:\RMod_{\DD_1(R)}\simeq \Ind(\textup{LCoh}(R))$
(cf. the proof of Lemma~\ref{leftcomplemmma}, \cite[4.8.5.8]{HA}).
By Lemma~\ref{dbar}, $\RMod_{\DD_1(R)}\simeq \RMod_{k\otimes_{\DD_2(R)}k}$.
Note that the $t$-structure on $\RMod_{\DD_1(R)}\simeq \RMod_{k\otimes_{\DD_2(R)}k}$
induced by the relative tensor product of $\Mod_k\otimes_{\RMod_{\DD_2(R)}}\Mod_k$ in $\PRTT$
is detemined by  the smallest full subcategory $\RMod_{\DD_1(R)}^{\le0}$
which contains $\DD_1(R)$ and is closed under small colimits and extensions. 
This full subcategory $\RMod_{\DD_1(R)}^{\le0}$
corresponds to $\Ind(\textup{LCoh}(R)^{\le0})=\Ind(\textup{LCoh}(R))^{\le0}$ through $\alpha$
(note that $\alpha$ carries $\DD_1(R)$ to $k$). Thus, $\alpha$ is an equivalence 
in $\PRTT$.
The argument similar to the proof of Lemma~\ref{leftcomplemmma} shows that
$\Ind(\textup{LCoh}(R))\to \LMod_R$ is a left completion functor (see \cite[Proposition 1.3.4]{Ind} for the case of commutative $R$).
Since $\RMod_{\DD_1(R)}\simeq \Ind(\textup{LCoh}(R))\to \LMod_R$ is left Kan extension of  $\Perf_{\DD_1(R)}\simeq \textup{LCoh}(R)$,
the composite $\RMod_{\DD_1(R)}\to \LMod_R$ is given by $N\otimes_{\DD_1(R)}k=(R\otimes_kN)\otimes_{(R\otimes_k\DD_1(R))}k$.

Next, we prove the latter assertion. 
To this end, we first note that
the left module action $\LMod_R^\otimes$ on $\RMod_{k\otimes_{\DD_2(R)}k}$ factors through 
the left module action of $\End_{\DD_2(R)}^l(\Mod_k)^\otimes$ on $\Mod_k\otimes_{\RMod_{\DD_2(R)}}\Mod_k\simeq \RMod_{k\otimes_{\DD_2(R)}k}$.
Here we use the notation in Lemma~\ref{dualcatalg}. By Lemma~\ref{dualcatalg},
there exists a $k$-linear monodial functor $\xi:\LMod_R^\otimes\simeq \LMod_{\DD_2(\DD_2(R))}^\otimes \to \End_{\DD_2(R)}^l(\Mod_k)^\otimes$.
The underlying functor $\LMod_R\to \End_{R}^l(\Mod_k)\simeq \RMod_{\DD_1(R)}$ carries 
$M$ to $k\otimes_{R}M\simeq k\otimes_{(R\otimes \DD_1(R))}(M\otimes \DD_1(R))$ (that is induced by the universal pairing $R\otimes_k\DD_1(R)\to k$).
(Since $R$ is Artin, this functor is fully faithful.
Indeed, $\xi|_{\Perf_R}:\Perf_R\to  \RMod_{\DD_1(R)}$ is fully faithful since
$R^{op}\simeq \End_{\LMod_R}(R)\to \End_{\RMod_{\DD_1(R)}}(k)\simeq \DD_1(\DD_1(R^{op}))$ is a canonical equivalence.
Moreover, $\xi(R)=k\in \RMod_{\DD_1(R)}$ is compact since $k\otimes_{\DD_1(R)}k\in \textup{LCoh}(R)$,
where we identify $k\otimes_{\DD_1(R)}k$ with the $k$-linear dual of $R\simeq \DD_1(\DD_1(R))$.
It follows that $\xi(\Perf_R)$ is the full subcategory of compact objects so that $\xi$ is fully faithful.)
We have the left module action functor $\LMod_R\times \RMod_{k\otimes_{\DD_2(R)}k}\to \RMod_{k\otimes_{\DD_2(R)}k}$.
To prove the assertion, it will suffice to prove that
there exists an object $E\in \RMod_{k\otimes_{\DD_2(R)}k}$ such that the composite 
\[
\rho:\LMod_R\times \{E\}\to \LMod_R\times \RMod_{k\otimes_{\DD_2(R)}k}\to  \RMod_{k\otimes_{\DD_2(R)}k}  \stackrel{\eta}{\to} \LMod_R
\]
is an equivalence.
We take $E$ to be the Hom complex $\Hom_R(k,R)\in \RMod_{k\otimes_{\DD_2(R)}k}\simeq \RMod_{\DD_1(R)}$, which corresponds to $R\in \Ind(\textup{LCoh}(R))$ lying over $R\in \LMod_R$ through $\eta$.
Since the (standard) free module $R\in \LMod_R$ is a single compact generator of $\LMod_R$, it is enough to show that $\rho$ is fully faithful when $\rho$ is restricted to the full subcategory spanned
by $R\in \LMod_R\simeq \LMod_R\times\{E\}$. 
The fully faithful functor $\LMod_R\to \RMod_{k\otimes_{\DD_2(R)}k}\simeq \RMod_{\DD_1(R)}$ carries $R$ to $k$ so that it
induces the equivalence of algebras $R^{op}\simeq \End_{\LMod_R}(R)\simeq \End_{\DD_1(R)^{op}}(k)=\End_{\RMod_{\DD_1(R)}}(k)$.
Here $\End_{\LMod_R}(R)$ and $\End_{\RMod_{\DD_1(R)}}(k)$ mean the endomorphism algebras defined as objects in $\Alg_1(\Mod_k)$.
The left action 
$\End_{\RMod_{\DD_1(R)}}(k)$ on $k$ induces the right action of $R$ on $R\simeq E\otimes_{\DD_1(R)}k$ which is the 
left action of $\End_{\LMod_R}(R)$ on $\rho(R)=R$.
By the definition of $E$, the standard left $R$-module structure 
on $R$ is (equivalent to) the left $R$-module structure $E\otimes_{\DD_1(R)}k\simeq (R\otimes_kE)\otimes_{(R\otimes\DD_1(R))}k\simeq R$ canonically determined by the universal pairing $R\otimes\DD_1(R)\to k$.
Since $R$ is an $\etwo$-algebra, $R\simeq R^{op}$ and $\DD_1(R)\simeq \DD_1(R)^{op}$ as $\eone$-algebras so that 
$R^{op}\otimes\DD_1(R)^{op}\simeq R\otimes\DD_1(R)\to k$. 
Thus, the left module action of $R^{op}\simeq \End_{\LMod_R}(R)$ on $\rho(R)=R$ is equivalent to the standard left module action of $R$ through $R\simeq R^{op}$.
Consequently, $\LMod_R\times\{E\}\to \LMod_R$ is an fully faithful when $\rho$ is restricted to the full subcategory spanned
by $(R, E)\in \LMod_R\times\{E\}$.
\QED


\subsection{}
Let $R$ be an augmented $\etwo$-algebra in $\Mod_k$, that is, an object of $\Alg_2^+(\Mod_k)$.
Suppose that $R$ is Artin.
Let $\lambda:R\otimes_k\DD_2(R)\to k$ be the Koszul dual pairing (cf. Section~\ref{koszulsec}).
Taking module categories, we have a monoidal functor  
\[
\LMod_R\otimes_k\LMod_{\DD_2(R)}\to \LMod_k=\Mod_k.
\]
Moreover, it preserves small colimits and carries $R\otimes_k\DD_2(R)$ to $k$.
Thus, it can be regarded as a morphism in $\Alg_1((\PRTT)_k)$.
In particular, $\Mod_k$ is a $\LMod_R^\otimes$-$\LMod_{\DD_2(R)}^{\otimes op}$ bimodule in $(\PRTT)_k$.
Recall that the superscript $\otimes op$ means the opposite algebra of $\LMod_{\DD_2(R)}^{\otimes}$
in $\Alg_1((\PRTT)_k)$.
This bimodule induces
\[
\xymatrix{
(-)\otimes_{R}k :\RMod_{\LMod_{R}}(\PRTT) \ar[r] &  \LMod_{\LMod_{\DD_2(R)}}(\PRTT) 
}
\]
where $(-)\otimes_{R}k$ indicates the functor induced by 
the relative tensor product over $\LMod_R$.
We write $\Phi_R$ for $(-)\otimes_{R}k$.
Let 
\[
\Psi_R :\LMod_{\LMod_{\DD_2(R)}}(\PRTT)  \to \RMod_{\LMod_{R}}(\PRTT)
\]
be the functor which sends $\EE$ to $\Mod_k\otimes_{\LMod_{\DD_2(R)}} \EE$, where $\Mod_k$ is regarded as a right $\LMod_R\otimes_k\LMod_{\DD_2(R)}$-module object.

\begin{Lemma}
\label{moritaleft}
The functor $\Psi_R$ is a left adjoint of $\Phi_R$. 
\end{Lemma}

\Proof
We define
\[
U_{\EE}:\EE\to (\LMod_k\otimes_{\LMod_{\DD_2(R)}}\EE)\otimes_{\LMod_R}\RMod_k \simeq (\RMod_k\otimes_{\LMod_R^{\otimes op}}\LMod_k)\otimes_{\LMod_{\DD_2(R)}}\EE
\]
to be the morphism obtained from $\xi:\LMod_{\DD_2(R)}\to \Mod_k\otimes_{\RMod_R}\Mod_k\simeq \RMod_k\otimes_{\LMod_R^{\otimes op}}\LMod_k$
(see Lemma~\ref{leftcomplemmma} and Remark~\ref{moduleproduct2})
by taking the tensor product over $\LMod_{\DD_2(R)}$.
We define
\[
V_{\DDD}:\LMod_k\otimes_{\LMod_{\DD_2(R)}}(\DDD\otimes_{\LMod_{R}}\RMod_k) \simeq  \DDD\otimes_{\LMod_{R}}(\RMod_k\otimes_{\LMod_{\DD_2(R)}^{\otimes op}}\LMod_k)\to \DDD
\]
to be the morphism obtained from $\eta:\RMod_k\otimes_{\LMod_{\DD_2(R)}^{\otimes op}}\LMod_k\simeq \Mod_k\otimes_{\RMod_{\DD_2(R)}}\Mod_k\to \LMod_R$ (see Lemma~\ref{leftcomp2})
by taking the tensor product over $\LMod_{R}$.
We set $\Phi=\Phi_R$ and $\Psi=\Psi_R$.
They determine the natural transformations $U:\textup{id}\to \Phi\circ \Psi$ and $V:\Psi \circ \Phi\to \textup{id}$.
We show that $U$ and $V$ are compatible up to homotopy (\cite[2.1]{RV}, \cite[6.2.1]{Kero}).
To this end, it suffices to prove:
\begin{enumerate}
\item The identity functor $\textup{id}_{\Psi}$ is equivalent to the composition
\[
\Psi=\Psi\circ \textup{id} \stackrel{\textup{id}_{\Psi}\circ U}{\longrightarrow}  \Psi\circ \Phi\circ \Psi \stackrel{V\circ \textup{id}_\Psi}{\longrightarrow} \textup{id} \circ \Psi
\]

\item
The identity functor $\textup{id}_{\Phi}$ is equivalent to the composition
\[
\Phi=\Phi\circ \textup{id} \stackrel{U\circ \textup{id}_{\Phi}}{\longrightarrow}  \Phi\circ \Psi\circ \Phi \stackrel{\textup{id}_\Phi \circ V}{\longrightarrow} \textup{id} \circ \Phi.
\]

\end{enumerate}
The proof of (2) is similar to (1) so that we will concentrate (1).
Unwinding the composition in (1) it is enough to prove the identity functor $\Mod_k\to \Mod_k$ is equivalent to the composition  
\begin{eqnarray*}
c:\Mod_k\simeq\Mod_k\otimes_{\LMod_{\DD_2(R)}}\LMod_{\DD_2(R)}&\stackrel{s}{\to}& \Mod_k\otimes_{\LMod_{\DD_2(R)}}\Mod_k\otimes_{\RMod_R}\Mod_k \\
&\stackrel{t}{\to}&  \RMod_R \otimes_{\RMod_R}\Mod_k \simeq \Mod_k
\end{eqnarray*}
Here the functor $s$ is induced by $\xi$, and $t$
is induced by $\eta:\Mod_k\otimes_{\LMod_{\DD_2(R)}}\Mod_k\simeq \Mod_k\otimes_{\RMod_{\DD_2(R)}}\Mod_k\to \RMod_R$.
After passing to the left completion, $s$ can be identified with the equivalence
\[
\Mod_k\otimes_{\HL(\LMod_{\DD_2(R)})}\HL(\LMod_{\DD_2(R)})\to \Mod_k\otimes_{\HL(\LMod_{\DD_2})}\Mod_k\otimes_{\RMod_R}\Mod_k
\]
induced by $\HL(\LMod_{\DD_2(R)})\simeq \Mod_k\otimes_{\RMod_R}\Mod_k$.
After passing to the left completion, $t$ can be identified with the equivalence
\[
\HL(\Mod_k\otimes_{\LMod_{\DD_2(R)}}\Mod_k)\otimes_{\RMod_R}\Mod_k \to  \RMod_R \otimes_{\RMod_R}\Mod_k
\]
induced by $\HL(\Mod_k\otimes_{\LMod_{\DD_2(R)}}\Mod_k)\simeq \RMod_R$ (see Lemma~\ref{leftcomp2}).
We also observe that the $k$-linear composite $c$ carries $k\in \Mod_k$ in the domain maps to $k\in \Mod_k$ in the target
(it is enough to verify this after the left completion).
Note that the $t$-structure on $\Mod_k$ is left complete.
We conclude that the composition $c$ is equivalent to the identity. 
\QED

\begin{Proposition}
\label{unit}
Let $\EE$ be a left $\LMod_{\DD_2(R)}$-module in $\PRTT$.
Let $U_{\EE}:\EE\to \Phi_R\circ \Psi_R(\EE)$ be the morphism in $\LMod_{\LMod_{\DD_2(R)}}(\PRTT)$
which is induced by the unit map of the adjunction.
Then $U$ is an equivalence after a left completion.
\end{Proposition}

\Proof
By Lemma~\ref{moritaleft}, 
\[
\Phi_R\circ \Psi_R(\EE)\simeq (\LMod_k\otimes_{\LMod_{\DD_2(R)}}\EE)\otimes_{\LMod_R}\RMod_k.
\]
We have equivalences
\[
(\LMod_k\otimes_{\LMod_{\DD_2(R)}}\EE)\otimes_{\LMod_R}\RMod_k\simeq \RMod_k\otimes_{\LMod_R^{\otimes op}}\LMod_k\otimes_{\LMod_{\DD_2(R)}}\EE
\]
and
\begin{eqnarray*}
  (\RMod_k\otimes_{\LMod_R^{\otimes op}}\LMod_k)\otimes_{\LMod_{\DD_2(R)}}\EE &\simeq& (\RMod_k\otimes_{\RMod_R}\LMod_k)\otimes_{\LMod_{\DD_2(R)}}\EE 
\end{eqnarray*}
(cf. Proposition~\ref{oppositemonoidal}).
Let $\End_{\RMod_R}^l(\Mod_k)^\otimes\in \Alg_1(\PR_k)$ be the endomorphism algebra of $\Mod_k\in \LMod_{\RMod_A}(\PR_k)$.
The right action of $\LMod_{\DD_2(R)}$ on $\RMod_k\otimes_{\RMod_R}\LMod_k$ in the last relative tensor product
is obtained from the universal right module action of $\End_{\RMod_A}^l(\Mod_k)^{\otimes op}$ on $\Mod_k=\LMod_k$
and the monoidal left completion functor
\[
\LMod_{\DD_2(R)}^\otimes\simeq \RMod_{\DD_2(R)}^{\otimes op}\to \End_{\RMod_R}^l(\Mod_k)^{\otimes op}
\]
(cf. the proof of Lemma~\ref{dualcatalg}).
See Proposition~\ref{oppositemonoidal} for the equivalence on the left-hand side.
The right-hand side is induced by $\RMod_{\DD_2(R)}^{\otimes}\to \End_{\RMod_R}^l(\Mod_k)^{\otimes}$
(the right version of $\xi$ in Corollary~\ref{monoidal}, see also  Lemma~\ref{leftcomplemmma})
so that we can identify $\End_{\RMod_R}^l(\Mod_k)^{\otimes op}$ with $\End_{R}^l(\Mod_k)^{\otimes}=\End_{\LMod_R}^l(\Mod_k)^{\otimes}$.

Note that the left module action of $\End_{\RMod_R}^l(\Mod_k)^\otimes$ exhibits
$\Mod_k\otimes_{\RMod_R}\Mod_k$ as a free module of rank one.
To see this, it is enough to prove that
there exists an object $m\in \Mod_k\otimes_{\RMod_R}\Mod_k$
such that the module action functor $\End_{\RMod_R}^l(\Mod_k)\times \Mod_k\otimes_{\RMod_R}\Mod_k\to \Mod_k\otimes_{\RMod_R}\Mod_k$
induces an equivalence $\End_{\RMod_R}^l(\Mod_k)\times \{m\}\to \Mod_k\otimes_{\RMod_R}\Mod_k$.
There exists an equivalence 
$\End_{\RMod_R}^l(\Mod_k)\simeq \Mod_k\otimes_{\RMod_R}\Mod_k$
and the module action functor is induced by the restriction functor
\begin{eqnarray*}
(\Mod_k\otimes_{\RMod_R}\Mod_k)\otimes_k (\Mod_k\otimes_{\RMod_R}\Mod_k) &\simeq& \RMod_{(k\otimes_Rk) \otimes_k (k\otimes_Rk)} \\
&\to&  \RMod_{(k\otimes_Rk)} \\
&\simeq&  \Mod_k\otimes_{\RMod_R}\Mod_k 
\end{eqnarray*}
along the comultiplication morphism $k\otimes_Rk \to  (k\otimes_Rk) \otimes_k (k\otimes_Rk)$ in $\Alg_1(\Mod_k)$
given by the bar construction of $R$.
When we regard $k$ as a $(k\otimes_Rk)$-module, the restriction
\[
(\Mod_k\otimes_{\RMod_R}\Mod_k)\simeq (\Mod_k\otimes_{\RMod_R}\Mod_k)\times \{k\}\to (\Mod_k\otimes_{\RMod_R}\Mod_k)
\]
is equivalent to the identity functor because the composition
$(k\otimes_Rk) \to  (k\otimes_Rk) \otimes_k (k\otimes_Rk)\to (k\otimes_Rk)\otimes_k k$
is the identity where the second morphism is induced by the counit $(k\otimes_Rk)\to k$.
We conclude that $\Mod_k\otimes_{\RMod_R}\Mod_k$ is a free $\End_{\RMod_A}^l(\Mod_k)^\otimes$-module 
of rank one.
The left completion functor $\HL$ is a symmetric monoidal functor which preserves small colimits so that we have
\[
\HL(\RMod_{k\otimes_Rk}\otimes_{\LMod_{\DD_2(R)}}\EE)\simeq \HL(\HL(\RMod_{k\otimes_Rk})\otimes_{\HL(\LMod_{\DD_2(R)})}\EE)\simeq \HL(\RMod_{k\otimes_Rk}\otimes_{\End_{R}^l(\Mod_k)}\EE)\simeq \EE.
\]
By the construction (cf. the proof of Lemma~\ref{moritaleft}), $U$ is induced by the left completion functor $\LMod_{\DD_2(R)}\to \Mod_{k}\otimes_{\RMod_k}\Mod_k$ regarded as the left $\LMod_{\DD_2(R)}$-module map
(namely, it is given by $\LMod_{\DD_2(R)}\otimes_{\LMod_{\DD_2(R)}}(-)\to (\Mod_{k}\otimes_{\RMod_k}\Mod_k)\otimes_{\LMod_{\DD_2(R)}}(-)$).
This proves our assertion.
\QED

Here is an application to the deformation theory of stable $\infty$-categories without $t$-structures.

\begin{Corollary}
\label{without}
Let $\CCC$ be a $k$-linear stable presentable $\infty$-category, that is, an object of $\PR_k$.
Suppose that $\CCC$ has a left and right complete, accessible $t$-structure $(\CCC^{\le0},\CCC^{\ge0})$.
Let $A$ be an augmented $\etwo$-algebra in $\Mod_k$.
Suppose that $\CCC$ is a promoted to $\CCC'\in \LMod_{\LMod_{\DD_2(A)}}(\PR)$ (this amounts to giving a map $f:\DD_2(R)\to \HH^\bullet(\CCC/k)$, see Proposition~\ref{tomapping}).
Then the left completion 
\[
\widehat{L}(\Mod_k\otimes_{\LMod_{\DD_2(A)}}\CCC')\in \RMod_{\LMod_{R}}(\PR)
\]
gives a deformation of $\CCC$ to $R$ (without $t$-structures).
Namely, there exists a canonical equivalence $\widehat{L}(\Mod_k\otimes_{\LMod_{\DD_2(A)}}\CCC')\otimes_{\LMod_R}\Mod_k\simeq \CCC$ in $\PR_k$.
If we regard  $\widehat{L}(\Mod_k\otimes_{\LMod_{\DD_2(A)}}\CCC')\otimes_{\LMod_R}\Mod_k$ 
as the induced left $\LMod_{\DD_2(A)}$-module, it is naturally equivalent to $\CCC'$ corresponding to $f$.
\end{Corollary}

\Proof
Use Lemma~\ref{already}, Proposition~\ref{unit}, and the natural transformation $U_{\CCC'}$.
\QED

\begin{Proposition}
\label{counit}
Let $\DDD$ be a left $\LMod_R$-module in $\PRTT$.
Let $V_{\DDD}:\Psi_R\circ \Phi_R(\DDD)\to \DDD$ be the morphism in $\RMod_{\LMod_{R}}(\PRTT)$
which is induced by the counit map of the adjunction.
Then $V_{\DDD}$ is an equivalence after a left completion.
\end{Proposition}

\Proof
By definition, 
\[
\Psi_R\circ \Phi_R(\DDD)\simeq \LMod_k\otimes_{\LMod_{\DD_2(R)}}(\DDD\otimes_{\LMod_{R}}\RMod_k).
\]
Moreover, there exist canonical equivalences
\begin{eqnarray*}
\LMod_k\otimes_{\LMod_{\DD_2(R)}}(\DDD\otimes_{\LMod_{R}}\RMod_k) &\simeq&  \DDD\otimes_{\LMod_{R}}(\RMod_k\otimes_{\LMod_{\DD_2(R)}^{\otimes op}}\LMod_k) \\
&\simeq& \DDD\otimes_{\LMod_{R}}(\RMod_k\otimes_{\RMod_{\DD_2(R)}}\LMod_k) \\
&\simeq& \DDD\otimes_{\LMod_{R}}  \RMod_{k\otimes_{\DD_2(R)}k}.
\end{eqnarray*}
By Lemma~\ref{leftcomp2}, there exists
a left completion functor $\RMod_{k\otimes_{\DD_2(R)}k}\to \LMod_R$.
It follows that
\begin{eqnarray*}
\HL(\DDD\otimes_{\LMod_{R}}  \RMod_{k\otimes_{\DD_2(R)}k}) &\simeq& \HL(\HL(\DDD)\otimes_{\HL(\LMod_{R})}\HL(\RMod_{k\otimes_{\DD_2(R)}k}))) \\
&\simeq& \HL(\HL(\DDD)\otimes_{\LMod_R}\LMod_R) \\
&\simeq& \HL(\DDD).
\end{eqnarray*}
By the construction (cf. Lemma~\ref{moritaleft}), the counit map is induced by $\RMod_{k\otimes_{\DD_2(R)}k}\to \LMod_R$.
Thus, our assertion follows.
\QED

\section{Main results}
\label{mainsec}

We prove main results of this paper.

\subsection{}
Let $(\EE,\EE^{\le0})\in (\PRTT)_k$ be a $k$-linear stable presentable $\infty$-category $\EE$
equipped with a right complete, accessible $t$-structure $(\EE^{\le0},\EE^{\ge0})$.
Let $\End_k(\EE)^\otimes\in \Alg_1(\PR_k)$ denote the endomorphism monoidal $\infty$-category
which has a universal left module action on $\EE$.

\begin{Definition}
\label{Hochdef}
We define the Hochschild cohomology $\etwo$-algebra $\HH^\bullet(\EE/k)$
to be $E_k(\End_k(\EE)^\otimes)\in \Alg_2(\Mod_k)$ (see Section~\ref{MOA} for $E_k:\Alg_1(\PR_k)\to \Alg_2(\Mod_k)$).
The underlying $\eone$-algebra is the endomorphism algebra of the identity functor in $\End_k(\EE)^\otimes$.
\end{Definition}

\begin{Remark}
One may consider an accessible $t$-structure on $\End_k(\EE)$ determined by the presentable full subcategory $\End_k(\EE)^{\le0}$
spanned by those $k$-linear functors $F$ such that $F(\EE^{\le0})\subset \EE^{\le0}$.
Namely,
$(\End_k(\EE),\End_k(\EE)^{\le0})$ is an endomorphism algebra object defined in $\Alg_1(\PRT)$.
Moreover, by \cite[VIII, Corollary 4.6.7]{DAG}, $\End_k(\EE)^{\le0}$ can be identified with the $\infty$-category of right $t$-exact functors $\EE\to \EE$.
The identity functor $\EE\to \EE$ lies in $\End_k(\EE)^{\le0}$. The endomorphism algebra of the identity functor
is the ordinary Hochschild cohomology $\HH^\bullet(\EE/k)$ (defined above),
and the canonical monoidal functor 
$\LMod_{\HH^\bullet(\EE/k)}^\otimes\to \End_k(\EE)^\otimes$
is right $t$-exact when $\LMod_{\HH^\bullet(\EE/k)}$ is equipped with the $t$-structure given in Remark~\ref{unifiedtstr}.
If one regards the Hochschild cohomology as the endomorphism algebra object of the identity,
it does not depend on the $t$-structure on $\EE$ (thus we do not use the $t$-structure on $\EE$ in Definition~\ref{Hochdef}).
\end{Remark}

Let 
$\LMod_{\LMod_{\DD_2(R)}}(\PRTT)^{\simeq}\times_{(\PRTT)_k^{\simeq}}\{\EE\}$ be the fiber product in $\wCat$
where $q:\LMod_{\LMod_{\DD_2(R)}}(\PRTT)\to (\PRTT)_k$ is the restriction functor 
along the monoidal functor $\Mod_k\to \LMod_{\DD_2(R)}$. As usual, by abuse of notation we write $\EE$ for $(\EE,\EE^{\le0})$.

\begin{Proposition}
\label{tomapping}
There exists a canonical equivalence
\[
\LMod_{\LMod_{\DD_2(R)}}(\PRTT)^{\simeq}\times_{(\PRTT)_k^{\simeq}}\{\EE\}\simeq \Map_{\Alg_2(\Mod_k)}(\DD_2(R), \HH^\bullet(\EE/k)).
\]
\end{Proposition}

\Proof
We first observe that any left $\LMod_{\DD_2(R)}$-module action on $\EE$ in $\PR_k$
is a left $\LMod_{\DD_2(R)}$-module action on $\EE$ in $(\PRTT)_k$ with respect to 
the $t$-structures on $\LMod_{\DD_2(R)}$ and $\EE$.
To see this, it is enough to show that 
the action functor $f:\LMod_{\DD_2(R)}\times \EE\to \EE$
sends $\LMod_{\DD_2(R)}^{\le0}\times \EE^{\le0}$
to $\EE^{\le0}$.
Let $\End_k(\EE)\in \Alg_1(\PR_k)$ be the endomorphism monoidal $\infty$-category
which has a universal left module action on $\EE$.
The left $\LMod_{\DD_2(R)}$-module action on $\EE$
factors through an essentially unique morphism $\LMod_{\DD_2(R)}\to \End_k(\EE)$ in $\Alg_1(\PR_k)$.
We note that $\LMod_{\DD_2(R)}\to \End_k(\EE)$ sends $\DD_2(R)$ to the identity functor of $\EE$.
Note that $f$ preserves small colimits in each variable, and
$\EE^{\le0}$ is closed under small colimits and extensions.
Taking into account  the definition of $\LMod_{\DD_2(R)}^{\le0}$, 
we see that $\LMod_{\DD_2(R)}^{\le0}\times \EE^{\le0}$
maps to $\EE^{\le0}$.
Consequently, 
we are reduced to showing that there exists a canonical equivalence
\[
\LMod_{\LMod_{\DD_2(R)}}(\PR_k)^{\simeq}\times_{(\PR_k)^{\simeq}}\{\EE\}\simeq \Map_{\Alg_2(\Mod_k)}(\DD_2(R), \HH^\bullet(\EE/k)).
\]
By the universality of $\End_k(\EE)$, we have an equivalence
\[
\LMod_{\LMod_{\DD_2(R)}}(\PR_k)^{\simeq}\times_{(\PR_k)^{\simeq}}\{\EE\}\simeq \Map_{\Alg_1(\PR_k)}(\LMod_{\DD_2(R)},\End_k(\EE)).
\]
Using the fully faithful functor $\Theta_k:\Alg_2(\Mod_k)\to \Alg_1(\PR_k)$ (that is, the left adjoit of $E_k$), we have the equivalence
\[
\Map_{\Alg_2(\Mod_k)}(\DD_2(R),\HH^\bullet(\EE/k))\simeq \Map_{\Alg_1(\PR_k)}(\LMod_{\DD_2(R)},\End_k(\EE)).
\]
This completes the proof.
\QED

We consider the sequence of the functors
\[
\RMod_{\LMod_{R}}(\PRTT) \stackrel{\Phi_{R}}{\longrightarrow}  \LMod_{\LMod_{\DD_2(R)}}(\PRTT) \stackrel{q}{\longrightarrow} (\PRTT)_k
\]
where $q$ is the forgetful functor (i.e. the restriction along $\Mod_k\to \LMod_{\DD_2(R)}$).
We shall write $\beta_{R\to k}^{+}:\RMod_{\LMod_{R}}(\PRTT) \to (\PRTT)_k$ for the base change functor along $\LMod_R\to \Mod_k$
in $\PRTT$, which is given by the base change in $\PR_{\mathbb{S}}$ at the level of underlying presentable $\infty$-categories.
Let $\beta_{R\to k}^{\pm}:\RMod_{\LMod_{R}}(\PRTTT) \to (\PRTTT)_k$ be the base change functor along $\LMod_R\to \Mod_k$
in $\PRTTT$. By Proposition~\ref{already}, it can be identified with the restriction of $\beta_{R\to k}^+$.  
We put
\[
\textup{Deform}_{\EE}(R)=\RMod_{\LMod_{R}}(\PRTTT)^\simeq \times_{\beta_{r\to k}^{+},(\PRTTT)^{\simeq}_k}\{\EE\}
\]
where the superscript $\simeq$ indicates the largest Kan subcomplex of the $\infty$-category: $\RMod_{\LMod_{R}}(\PRTTT)^\simeq$ is the largest $\infty$-groupoid contained in the whole $\infty$-category.
An object of $\textup{Deform}_{\EE}(R)$ can be regarded as a pair $(\DDD, \beta_{R\to k}^+(\DDD)\simeq \EE)$
such that $\DDD$ is a complete $R$-linear stable $\infty$-category, and an equivalence
$\DDD\otimes_{\LMod_R}\Mod_k=\beta_{R\to k}^+(\DDD)\simeq \EE$ in $\PRTT$.
We shall refer to  an object of $\textup{Deform}_{\EE}(R)$
is a complete deformation of $\EE$ to $R$.

\begin{Example}
Let
$X$ be a quasi-compact and separated derived (spectral) scheme over $k$. Let $R\to k$ be an argumented commutative connective algebra in $\Mod_k$.
Suppose that $X'\to \Spec R$ is a deformation of $X$ to $R$: There exists an equivalence $f:X'\times_{\Spec R}\Spec k\simeq X$ over $k$.
For a derived scheme $Y$, we let $(\textup{QCoh}(Y), \textup{QCoh}(Y))^{\le0})$ denote the stable $\infty$-category of quasi-coherent sheaves on $Y$ equipped with the left and right complete, accessible $t$-structure (see \cite[VIII]{DAG}). 
Then 
\[
\big{(}(\textup{QCoh}(X'), \textup{QCoh}(X')^{\le0}), f^*:\textup{QCoh}(X')\otimes_{\LMod_R}\Mod_k\simeq \textup{QCoh}(X)\big{)}
\]
provides a complete deformation of $(\textup{QCoh}(X),\textup{QCoh}(X)^{\le0})$ to $R$.
\end{Example}

Let $\EE\in (\PRTTT)_k$.
By Proposition~\ref{already}, $\Phi_R$ induces the functor 
\[
\rho:\RMod_{\LMod_{R}}(\PRTTT)\times_{\beta_{R\to k}^+,(\PRTT)_k}\{\EE\} \longrightarrow \LMod_{\LMod_{\DD_2(R)}}(\PRTTT)\times_{q,(\PRTTT)_k}\{\EE\}.
\]
According to Proposition~\ref{tomapping}, $\rho$ induces
\[
P_R:\textup{Deform}_{\EE}(R)\to \Map_{\Alg_2(\Mod_k)}(\DD_2(R), \HH^\bullet(\EE/k)).
\]

\begin{Theorem}
\label{mainthm}
The functor $P_R:\textup{Deform}_{\EE}(R)\to \Map_{\Alg_2(\Mod_k)}(\DD_2(R), \HH^\bullet(\EE/k))$
is an equivalence of spaces.
\end{Theorem}

\begin{Example}
\label{exsqzero}
Let $k\oplus k[n]\in \Alg_2^+(\Mod_k)$ denote the trivial square zero extension of $k$ by $k[n]$, which is endowed with the augmentation $k\oplus k[n]\to k$ (cf. \cite[X]{DAG}).
Then $\DD_2(k\oplus k[n])\in \Alg_2(\Mod_k)$ is a free $\etwo$-algebra generated by $k[-n-2]$ (see \cite[X, Proposition 4.5.6]{DAG}).
Thus, $\pi_0(\Map_{\Alg_2(\Mod_k)}(\DD_2(k\oplus k), \HH^\bullet(\EE/k))\simeq \HH^{n+2}(\EE/k))$.
The set $\pi_0(\textup{Deform}_{\EE}(k\oplus k[n]))$ of isomorphism classes of deformations to $k\oplus k[n]$ is naturally isomorphic to $\HH^{n+2}(\EE/k)=H^{n+2}(\HH^\bullet(\EE/k))$.
\end{Example}

{\it Proof of Theorem~\ref{mainthm}.}
The symmetric monoidal left completion functor $\HL:\PRTT\to \PRTTT$
induces
$\RMod_{\LMod_R}(\PRTT)\to \RMod_{\LMod_R}(\PRTTT)$
since the $t$-structure on $\LMod_R$ is left complete.
Consider 
\[
\RMod_{\LMod_R}(\PRTTT)\to \RMod_{\LMod_R}(\PRTT) \stackrel{\beta_{R\to k}^+}{\longrightarrow} (\PRTT)_k
\]
where the functor on the left-hand side is induced by the lax symmetric monoidal inclusion
$\PRTTT\to \PRTT$, and the right-hand side is defined by $\beta_{R\to k}^+=\otimes_{\LMod_R}\Mod_k$.
We abuse notation by writing $\beta_{R\to k}^+$ for the composite.
By Proposition~\ref{already}, this composite $\beta_{R\to k}^+$ is equivalent to
\[
\RMod_{\LMod_R}(\PRTTT)\stackrel{\beta_{R\to k}^{\pm}}{\longrightarrow} (\PRTTT)_k\hookrightarrow (\PRTT)_k.
\]
Thus, we have
\[
\xymatrix{
\RMod_{\LMod_R}(\PRTT)\ar[r]^{\HL_R} \ar[d]_{\beta^+_{R\to k}} &  \RMod_{\LMod_R}(\PRTTT) \ar[d]^{\beta^{+}_{R\to k}} \\ 
(\PRTT)_k \ar[r]^{\HL_k}  &  (\PRTTT)_k 
}
\]
which commutes up to canonical homotopy,
where
horizontal arrows are induced by completion functors.
We consider the composition
\[
 \HL_k \circ \beta_{R\to k}^+ \circ \Psi_R  :\LMod_{\LMod_{\DD_2(R)}}(\PRTT)  \to (\PRTTT)_k.
\]
The composite
$\beta_{R\to k}^+ \circ \Psi_R$
is equivalent to $\Phi_R\circ \Psi_R$ after composing the forgetful functor $q$.
Thus, Proposition~\ref{unit} shows that  
$\HL_k \circ \beta_{R\to k}^+ \circ \Psi_R$
and $\beta_{R\to k}^+\circ  \HL_R\circ \Psi_R$
are equivalent to $\HL_k\circ q$.
We then obtain the diagram
\[
\xymatrix{
\RMod_{\LMod_R}(\PRTTT) \ar[r]^{\Phi_R} \ar[d]^{\beta_{R\to k}^+} & \LMod_{\LMod_{\DD_2(R)}}(\PRTTT) \ar[r]^{\HL_R \circ \Psi_R} \ar[d]^{ q}  &  \RMod_{\LMod_R}(\PRTTT) \ar[d]^{\beta^+_{R\to k}}   \\
(\PRTTT)_k  \ar[r]^{\textup{id}} &    (\PRTTT)_k  \ar[r]^{\textup{id}} &  (\PRTTT)_k
}
\]
which commutes up to homotopy. 
Here $\HL_R \circ \Psi_R$ in the diagram means the restricted functor $\HL_R \circ \Psi_R|_{\LMod_{\LMod_{\DD_2(R)}}(\PRTTT)}$.
By Proposition~\ref{counit}, the upper horizontal composite $\HL_R \circ \Psi_R\circ \Phi_R$ is naturally equivalent to the identity functor.
It gives rise to the sequence
\begin{eqnarray*}
\RMod_{\LMod_R}(\PRTTT)\times_{\beta_{R\to k}^+,(\PRTTT)_k}\{\EE\} &\stackrel{\rho}{\longrightarrow}& \LMod_{\LMod_{\DD_2(R)}}(\PRTTT)\times_{q,(\PRTTT)_k}\{\EE\} \\
&\longrightarrow& \RMod_{\LMod_R}(\PRTTT) \times_{(\PRTTT)_k}\{\EE\}
\end{eqnarray*}
whose composition is naturally equivalent to the identity functor,
where the second functor is 
induced by
$\beta_{R\to k}^+\circ  \HL_R\circ \Psi_R$ and $\HL_R\circ \Psi_R$.

Next, we consider 
\begin{eqnarray*}
\LMod_{\LMod_{\DD_2(R)}}(\PRTTT)\times_{q,(\PRTTT)_k}\{\EE\} &\longrightarrow& \RMod_{\LMod_R}(\PRTTT) \times_{\beta_{R\to k}^+,(\PRTTT)_k}\{\EE\} \\
&\stackrel{\rho}{\longrightarrow}& \LMod_{\LMod_{\DD_2(R)}}(\PRTTT)\times_{q,(\PRTTT)_k}\{\EE\}.
\end{eqnarray*}
It will suffice to prove that the composite is naturally equivalent to the identity functor.
It is enough to show that $\Phi_R\circ \HL_R \circ \Psi_R:\LMod_{\LMod_{\DD_2(R)}}(\PRTT)\to \LMod_{\LMod_{\DD_2(R)}}(\PRTTT)$ is naturally equivalent to the functor $\HL_{\LMod_{\DD_2(R)}}:\LMod_{\LMod_{\DD_2(R)}}(\PRTT)\to \LMod_{\LMod_{\DD_2(R)}}(\PRTTT)$ induced by 
$\HL_k$ and the restriction along $\LMod_{\DD_2(R)}\to \HL_k(\LMod_{\DD_2(R)})$.
According to Proposition~\ref{already} and the fact that $\HL$ is a symmetric monoidal colimit-preserving functor, we see that the diagram
\[
\xymatrix{
\RMod_{\LMod_R}(\PRTT) \ar[r]^{\HL_R} \ar[d]^{\Phi_R}  &  \RMod_{\LMod_R}(\PRTTT) \ar[d]^{\Phi_R} \\
\LMod_{\LMod_{\DD_2(R)}}(\PRTT) \ar[r]_{\HL_{\LMod_{\DD_2(R)}}} & \LMod_{\LMod_{\DD_2(R)}}(\PRTTT)
}
\]
commutes up to canonical homotopy. It follows from Proposition~\ref{unit} that $\Phi_R\circ \HL_R \circ \Psi_R\simeq \HL_k\circ \Phi_R\circ \Psi_R\simeq \HL_{\LMod_{\DD_2(R)}}$.
This completes the proof.
\QED

\begin{Remark}
For a morphism $f:\DD_2(R)\to \HH^{\bullet}(\EE/k)$ of $\etwo$-algebras,
the corresponding deformation is defined as $\widehat{L}(\Mod_k\otimes_{\LMod_{\DD_2(R)}}\EE)$.
Here the left module action of $\LMod_{\DD_2(R)}^\otimes$ on $\EE$ is induced by $f$, the canonical monoidal functor
$\LMod_{\HH^{\bullet}(\EE/k)}^\otimes\to \mathcal{E}nd_{k}(\EE)^\otimes$, and the canonical action of $\mathcal{E}nd_{k}(\EE)^\otimes$ on $\EE$.
\end{Remark}

\subsection{}
We interpret Theorem~\ref{mainthm} in terms of formal moduli problems, that is, pointed formal stacks.
We first define a functor which describes the deformation problem of $\EE\in (\PRTTT)_k$.

\begin{Definition}
\label{deffunctorcomplete}
Let 
$\Theta_k^\pm:\Alg_2(\Mod_k^{\le0})\to \Alg_1((\PRTTT)_k)$ be the functor given by $B\mapsto (\LMod_B^\otimes,\LMod_B^{\le0})$ (see Construction~\ref{thetaTTT}). We set
\[
\textup{RLin}^{\wedge}(k):=\RMod((\PRTTT)_k)\times_{\Alg_1((\PRTTT)_k),\Theta_k^\pm}\Alg_2(\Mod_k^{\le0}).
\]
where $\RMod((\PRTTT)_k)\to \Alg_1((\PRTTT)_k)$ is the forgetful functor.
Let
 $\textup{RLin}^\wedge(k)^\dagger$ be the subcategory spanned by coCartesian morphisms over $\Alg_2(\Mod_k^{\le0})$.
Let
$\textup{RLin}^\wedge(k)_{/(\EE,k)}^\dagger$
be the overcategory where $(\EE,k)$ means $\EE\in (\PRTTT)_k$ lying over $k\in \Alg_2(\Mod_k^{\le0})$.
We then have a left fibration $r:\textup{RLin}^\wedge(k)_{/(\EE,k)}^\dagger\to \Alg_2^+(\Mod_k^{\le0})$.
Given $R\in \Alg_2^+(\Mod_k^{\le0})$ we refer to an object of the fiber $r^{-1}(R)$ over $R$ as a complete deformation of $\EE$ to $R$. 

We define 
\[
\overline{\textup{Deform}}_{\EE}:\Alg_2^+(\Mod_k^{\le0})\to \wSSS
\]
to be the functor which corresponds to the left fibration $r$. 
Let $\textup{Deform}_{\EE}:\textup{Art}_2\to \wSSS$ be the restriction of $\overline{\textup{Deform}}_{\EE}$ to $\textup{Art}_2$.
\end{Definition}

Next, we recall formal $\eenu$-moduli problems from \cite[Section 4]{DAG}.
Let $n$ be a natural number or $\infty$.
A functor $F:\textup{Art}_n\to \SSS$ from $\infty$-category of Artin $\eenu$-algebras is a formal $\eenu$-moduli problem if
\begin{itemize}
  \item $F(k)$ is a contractible space,
\item for any pullback square \[
\xymatrix{
R \ar[r] \ar[d] &  R_0 \ar[d] \\
R_1 \ar[r] & R_{01}
}
\]
in $\textup{Art}_n$
such that $H^0(R_0)\to H^0(R_{01})$ and $H^0(R_1)\to H^0(R_{01})$ are surjective, the square
 \[
\xymatrix{
F(R) \ar[r] \ar[d] &  F(R_0) \ar[d] \\
F(R_1) \ar[r] & F(R_{01})
}
\]
is a pullback diagram.
\end{itemize}
Let $\widehat{\mathsf{Stack}}_n^\ast$ be the full subcategory of $\Fun(\textup{Art}_n, \SSS)$, which consists of formal $\eenu$-moduli problems.
When $n$ is a natural number, thanks to a theorem of Lurie \cite[X, Theorem 4.0.8]{DAG} there exists an equivalence of $\infty$-categories
\[
\widehat{\mathsf{Stack}}_n^\ast \simeq \Alg_n^+(\Mod_k).
\]
Through the equivalence,
an augmented algebra $B\in \Alg_n^+(\Mod_k)$ corresponds to the formal $\eenu$-moduli problem $\mathcal{F}_{B}^{(n)}: \textup{Art}_n\to  \SSS$
defined informally by $R \mapsto \Map_{\Alg_n^+(\Mod_k)}(\DD_n(R),B)$.

\begin{Theorem}
\label{mainthm2}
The deformation functor $\textup{Deform}_{\EE}:\textup{Art}_2\to \wSSS$
is a formal $\etwo$-moduli problem $\mathcal{F}_{k\oplus \HH^\bullet(\EE/k)}^{(2)}$ corresponding to the augmented $\etwo$-algebra $k\oplus \HH^\bullet(\EE/k)\to k$.
\end{Theorem}

\Proof
Consider the adjoint pair $U:\Alg_2^+(\Mod_k)\rightleftarrows \Alg_2(\Mod_k):T$ where $U$ is the forgetful functor, and $T$ is its right adjoint given by $A\mapsto [\textup{pr}_1:k\oplus A\to k]$
where $k\oplus A$ means the product of $\etwo$-algebras.
This adjoint pair induces a canonical equivalence of mapping spaces
$\eta_R:\Map_{\Alg_2(\Mod_k)}(\DD_2(R), \HH^\bullet(\EE/k))\simeq \Map_{\Alg_2^+(\Mod_k)}(\DD_2(R), k\oplus \HH^\bullet(\EE/k))$.
In view of Theorem~\ref{mainthm}, it is enough to construct 
the functor
\[
\textup{Deform}_{\EE}(R)\to  \Map_{\Alg_2^+(\Mod_k)}(\DD_2(R), k\oplus \HH^\bullet(\EE/k))
\]
functorially on $R\in \Alg_2^+(\Mod_k^{\le0})$, which is equivalent to $P_R$ through $\eta_R$.
Namely, it will suffice to construct a natural transformation $\textup{Deform}_{\EE}\to \mathcal{F}_{k\oplus \HH^\bullet(\EE/k)}^{(2)}$ having the required property.
The construction is done in \cite[X, Construction 5.3.18]{DAG} up to a small modification: One uses $(\PRTT)_k$ instead of  $\textup{LinCat}_k$ in {\it loc.cit.}
(strictly speaking, the convention on left and right modules is reversed with ours). 
Besides, we can also apply the axiomatic formulation in \cite[Section 4.4]{IM}.
\QED

\subsection{Commutative bases}
We consider the restriction to deformations to commutative bases. We assume that the base field $k$ is of characteristic zero.
Let $\CAlg(\Mod_k)\to \Alg_2(\Mod_k)$ be the canonical forgetful functor obtained from
the canonical map of $\infty$-operads $\mathbf{E}_2^\otimes\to \mathbf{E}_\infty^\otimes$.
It induces $\textup{Art}=\textup{Art}_{\infty}\to \textup{Art}_2$ so that the composition gives the functor
\[
\Fun(\textup{Art}_2,\SSS)\longrightarrow  \Fun(\textup{Art},\SSS).
\]
According to \cite[Lemma 7.5]{IM} it preserves formal moduli problems so that we have the restriction $res_{2/\infty}:\widehat{\mathsf{Stack}}_2^*\rightarrow \widehat{\mathsf{Stack}}_\infty^*$.
Let $\textup{Lie}_k$ be the $\infty$-category of dg Lie algebras over $k$.
The $\infty$-category is obtained from the category of dg Lie algebras by inverting quasi-isomorphisms.
Through the categorical equivalence $\textup{Lie}_k\simeq \widehat{\mathsf{Stack}}_\infty^*$ (see \cite[X, 2.0.2]{DAG}) and $\Alg_2^+(\Mod_k)\simeq \widehat{\mathsf{Stack}}_\infty^*$
the functor $res_{s/\infty}$ corresponds to 
\[
\Alg_2^+(\Mod_k)\to \textup{Lie}_k.
\]
The underlying complex of  $res_{2/\infty}(\epsilon:B\to k)$ is $\Ker(\epsilon)[1]$ where $\Ker$ denotes the fiber (see \cite[Proposition 7.7]{IM}).
By abuse of notation we write $\HH^\bullet(\EE/k)[1]$ for the underlying dg Lie algebra $res_{2/\infty}(k\oplus \HH^\bullet(\EE/k) \to k)$.
We set $\textup{Deform}_{\EE}^{\infty}=res_{2/\infty}(\textup{Deform}_{\EE})$.

\begin{Corollary} 
\label{mainthm3}
The deformation functor
$\textup{Deform}_{\EE}^{\infty}$ is a formal $\mathbf{E}_\infty$-moduli problem $\mathcal{F}_{\HH^\bullet(\EE/k)[1]}^{(\infty)}$ corresponding to $\HH^\bullet(\EE/k)[1]$.
\end{Corollary}

\subsection{Obstruction theory}
From a formal moduli problem, one can obtain the obstruction theory of deformations. We will explain this consequence for the reader's convenience.
We refer the reader to \cite[Section 3.3]{BT} for a comprehensive account.
Let $p:R'\to R$ be a morphism in $\textup{Art}_2$, which fits in the pullback square
\[
\xymatrix{
R' \ar[r]^p \ar[d] & R \ar[d] \\
k \ar[r] & k\oplus V[i] 
}
\]
in $\textup{Art}_2$, where $V$ is a finite dimensional vector space, $i$ is a positive integer,
and $k\oplus V[i]$ is the trivial square zero extension by $V[i]$

\begin{Example}
Let $p:R'\to R$ be a surjective map of ordinary (discrete) Artin local algebras in $\textup{Art}$.
Let $I=\Ker(p)$ and let $\mathfrak{m}_{R'}$ be the maximal ideal of $R'$.
We assume that  $\mathfrak{m}_{R'}\cdot I=0$, that is, $p$ is an elementary extension of Artin local algebras.
There exists a pullback square
\[
\xymatrix{
R' \ar[r]^p \ar[d] & R \ar[d] \\
k \ar[r] & k\oplus I[1] 
}
\]
in $\textup{Art}$, where $k\oplus I[1]$ is the trivial square zero extension by $I[1]$ (see e.g. \cite[Lemma 3.3]{BT}). 
We remark that every surjective map of ordinary Artin local algebras in $\textup{Art}$
decomposes as a sequence of elementary extensions.
\end{Example}

Suppose that  $\DDD\in \pi_0(\textup{Deform}_{\EE}(R))$.
We consider the problem of whether or not $\DDD$ is promoted to $R'$ along $p:R'\to R$.

\begin{Theorem}
\label{mainthm4}
There exists an element $\textup{Obs}(\DDD,p)\in \HH^{i+2}(\EE/k)\otimes_kV$ 
such that $\textup{Obs}(\DDD,p)=0$ if and only if there exists $\DDD'\in \pi_0(\textup{Deform}_{\EE}(R'))$ lying over $\DDD$ along  the map $p$.
\end{Theorem}

\Proof
By Theorem~\ref{mainthm2}, we have the canonical equivalence 
\[
\alpha:\textup{Deform}_{\EE}(R')\simeq \textup{Deform}_{\EE}(R)\times_{\textup{Deform}_{\EE}(k\oplus V[i])}\textup{Deform}_{\EE}(k).
\]
Note that $\textup{Deform}_{\EE}(k)$ is a contractible space.
Since $\DD_2(k\oplus V[i])$ is a free $\etwo$-algebra generated by $V^{\vee}[-i-2]$ (\cite[X, Proposition 4.5.6]{DAG}),
it follows that $\phi:\textup{Deform}_{\EE}(k\oplus V[i])\simeq \HH^{i+2}(\EE/k)\otimes_kV$.
Here, $V^\vee$ is the dual vector space.
Moreover, by the construction,
$0\in \HH^{i+2}(\EE/k)\otimes_kV$ corresponds to the trivial deformation, that is, the image of $\textup{Deform}_{\EE}(k)\to \textup{Deform}_{\EE}(k\oplus V[i])$.
We define $\textup{Obs}(\DDD,p)$ to be the image of $\DDD\in \pi_0(\textup{Deform}_{\EE}(R))$
in $\HH^{i+2}(\EE/k)\otimes_kV$ (through $\phi$).
Taking into account the equivalence $\alpha$, we see that 
there exists $\DDD'\in \pi_0(\textup{Deform}_{\EE}(R'))$ lying over $\DDD$ if and only if
$\textup{Obs}(\DDD,p)$ is equivalent to the trivial deformation.
The latter condition amounts to the condition that $\textup{Obs}(\DDD,p)=0$.
\QED

\end{document}